\documentclass[10pt]{article}

\usepackage[a4paper]{geometry}
\usepackage{graphicx}
\usepackage{caption}
\usepackage{indentfirst}
\usepackage[pdfpagemode=UseNone]{hyperref}
\usepackage{amssymb, amsmath, amsthm}

\newcommand{\di}{\mathbin{\lozenge}} 
\newcommand{\zp}{\mathord{\scalebox{0.675}[1]{$\mathbb{O}$}}} 

\DeclareMathOperator{\res}{Res}
\DeclareMathOperator{\disc}{Disc}

\newcommand{\mynewtheorem}[2]{\newtheorem{#1}{\indent #2}}
\newcommand{\myalttheorem}[2]{\newtheorem*{#1}{\indent #2}}
\mynewtheorem{conjecture}{Conjecture}
\myalttheorem{conjecture*}{Conjecture}
\mynewtheorem{lemma}{Lemma}
\mynewtheorem{proposition}{Proposition}
\myalttheorem{proposition*}{Proposition}
\mynewtheorem{theorem}{Theorem}
\myalttheorem{theorem*}{Theorem}
\newenvironment{myproof}[1][Proof]{\begin{proof}[\indent #1]}{\end{proof}}

\begin{document}

\title{\textbf{Diamond Determinants and Somos Sequences}}
\author{Nikolai Beluhov}
\date{}

\maketitle

\begin{abstract} A Somos sequence of order $n$ is defined by a quadratic recurrence of width $n + 1$. Some of the remarkable properties of these sequences for small $n$ are tied to certain matrices built out of them being of finite rank. We give an elementary proof of the finite-rank property for order $6$, previously only established with the help of advanced machinery from the theory of hyperelliptic functions. Our method also yields a new finite-rank property for the Somos sequences of order $7$. In addition, we conjecture generalisations of these results to higher orders, for the subclass of Gale-Robinson sequences. \end{abstract}

\section{Introduction} \label{intro}

Let $s$ be a sequence of complex numbers with $i$-th term $s_i$. Fix a positive integer $n$. Consider any window $s_i$, $s_{i + 1}$, $\ldots$, $s_{i + n}$ of size $n + 1$ in $s$, and form the products $s_{i + j}s_{i + n - j}$, with $j = 0$, $1$, $\ldots$, $\lfloor n/2 \rfloor$, of the symmetric pairs of elements in this window. We are interested in sequences $s$ where these products satisfy a fixed linear relation as the window slides back and forth along the sequence.

Formally, $s$ is a \emph{Somos sequence} of order $n$ when there exist constants $a_0$, $a_1$, $\ldots$, $a_{\lfloor n/2 \rfloor}$, with $a_0 \neq 0$, such that \[a_0s_is_{i + n} + a_1s_{i + 1}s_{i + n - 1} + \cdots + a_{\lfloor n/2 \rfloor}s_{i + \lfloor n/2 \rfloor}s_{i + \lceil n/2 \rceil} = 0\] for all $i$. Since scaling all of the $a$'s by the same nonzero factor simultaneously does not meaningfully change anything, from now on we are going to assume without loss of generality that $a_0 = -1$. Then \[s_{i + n} = \frac{a_1s_{i + 1}s_{i + n - 1} + \cdots + a_{\lfloor n/2 \rfloor}s_{i + \lfloor n/2 \rfloor}s_{i + \lceil n/2 \rceil}}{s_i},\] provided that $s_i \neq 0$; and we can similarly express $s_i$ in terms of $s_{i + 1}$, $s_{i + 2}$, $\ldots$, $s_{i + n}$ when $s_{i + n} \neq 0$. We care most of all about Somos sequences where all terms are nonzero, so that both the forward and the backward forms of this recurrence hold universally throughout the sequence.

Somos sequences exhibit some remarkable properties. Most of the literature focuses on orders $4$, $5$, $6$, $7$ as the properties in question are often true but trivial for lower $n$ and false for higher ones. The engagingly-written \cite{G} covers the early history of the subject, while the problem collection \cite{TT} offers an elementary survey of later developments. From a more advanced point of view, the study of Somos sequences has involved also cluster algebras \cite{FZ}, elliptic functions \cite{H05, H07}, and hyperelliptic functions \cite{H09, FH}. The parts of this background relevant to our work will be reviewed below.

We go on to consider one motivating problem. Let $d$ be a positive integer. When we take each $d$-th term of $s$, we obtain a subsequence of $s$ known as a \emph{decimation} of $s$ by a factor of~$d$. So, in particular, $s$ splits into $d$ such decimations. Surprisingly, the decimations of low-order Somos sequences turn out to be Somos sequences themselves -- though sometimes of an order higher than that of the original sequence. This behaviour is easy to observe experimentally, as we will see in Section \ref{init}. For the sake of clarity, initially we are going to limit ourselves to orders $6$ and $7$, with discussion of the lower orders $2$, $3$, $4$, $5$ postponed until Section \ref{low}.

We begin with order $6$. In this setting, some decimation properties (of a concrete Somos sequence, with decimation factors $6$ and $12$) were conjectured by Speyer \cite{S} regarding one number-theoretic problem. A proof was given by Ustinov \cite{U17} based on earlier work by Hone \cite{H09} as well as Fedorov and Hone \cite{FH} on the connections between Somos sequences and hyperelliptic functions. The argument shows, in essence, the following:

\begin{proposition} \label{o6d} Let $s$ be a Somos sequence of order $6$ with nonzero terms. Then every decimation of $s$ is a Somos sequence itself, of order at most $8$. \end{proposition} 

The same argument in fact implies one much stronger result of which Proposition \ref{o6d} is a quick corollary. Before we can state this result, though, we must set up some vocabulary.

Let $M$ be a matrix of complex numbers. Choose $r$ rows and $r$ columns of $M$. The entries of $M$ where these rows and columns meet form an $r \times r$ \emph{sub-matrix} of $M$ whose determinant is an $r \times r$ \emph{minor} of $M$. Suppose now that, instead of $r$ rows and $r$ columns, we choose $r$ diagonals and $r$ anti-diagonals which meet pairwise at $r^2$ entries of $M$. We say that these entries form a \emph{diamond sub-matrix} of $M$, and we call its determinant a \emph{diamond minor} of $M$.

The \emph{rank} of $M$ can be defined in two equivalent ways as the dimension of the linear hull of its rows or as the dimension of the linear hull of its columns. It is also well-known to equal the smallest nonnegative integer $r$ such that all minors in $M$ of size $(r + 1) \times (r + 1)$ vanish. By analogy, we define the \emph{diamond rank} of $M$ to be the smallest nonnegative integer $r$ such that all diamond minors in $M$ of size $(r + 1) \times (r + 1)$ vanish.

Let $t$ be a sequence of complex numbers with $j$-th term $t_j$. We write $s \times t$ for the matrix whose entry at position $(i, j)$ equals $s_it_j$. Suppose, temporarily, that both of $s$ and $t$ are doubly infinite. Then we can form the infinite matrices $M'$ and $M''$ whose entries at position $(i, j)$ are given by $s_{i - j}t_{i + j}$ and $s_{i - j}t_{i + j + 1}$, respectively. Clearly, the diamond sub-matrices of $s \times t$ coincide with the ordinary sub-matrices of $M'$ and $M''$. So the diamond rank of $s \times t$ equals the maximum of the ordinary ranks of $M'$ and $M''$.

Avdeeva and Bykovskii \cite{AB} define two doubly infinite sequences $s$ and $t$ to form a \emph{hyperelliptic system} of rank $r$ when this maximum equals $r$. The choice of term comes from the fact that such pairs of sequences arise naturally when we attempt to discretise certain addition formulas involving hyperelliptic functions. Beyond \cite{U17}, the applications of this framework to Somos sequences have been explored by Ustinov also in \cite{U19} and \cite{U24}.

We prefer the vocabulary of diamond sub-matrices and diamond ranks instead because it is better suited to our particular purposes; among other things, it will allow us to treat finite and infinite sequences in a uniform manner. We can now state the aforementioned stronger result:

\begin{theorem} \label{o6r} Let $s$ be a Somos sequence of order $6$ with nonzero terms. Then the matrix $s \times s$ is of diamond rank at most $4$. \end{theorem} 

Notice, though, that Theorem \ref{o6r} is a purely elementary statement. So it is somewhat odd that its only known proof should be based on such advanced machinery. One of our two main goals in the present paper will be to give a new proof of Theorem \ref{o6r} which does not venture outside of elementary linear algebra. The other one will be to employ our method so as to establish a similar finite-rank property also for the Somos sequences of order $7$.

We continue with a brief overview of this method. Suppose we wish to show that the rank of $M$ does not exceed $r$. One way to go about this would be to examine the $(r + 1) \times (r + 1)$ minors of $M$, and to verify that all of them vanish. Generally speaking, we cannot afford to miss even a single minor; it is possible to construct matrices where all but one $(r + 1) \times (r + 1)$ minors vanish, and yet the rank strictly exceeds $r$. However, under favourable circumstances -- provided that certain non-degeneracy conditions are satisfied -- we can get away with examining just the contiguous $(r + 1) \times (r + 1)$ minors of $M$. We spell out the details in Section \ref{cmin}.

So, in the setting of Theorem \ref{o6r}, we can afford to focus solely on the contiguous diamond minors of size $5 \times 5$ in $s \times s$. Once these have been tackled, the rest of the proof will amount to a bit of ``technical fiddling'' to ensure that the relevant non-degeneracy conditions are indeed satisfied.

Consider, then, any contiguous diamond minor of size $5 \times 5$ in $s \times s$. It is built out of the elements of two separate windows $s_X$ and $s_Y$ of size $9$ in $s$. Imagine, for a moment, some kind of hypothetical calculation which shows that this minor vanishes. How does this calculation ``know'' that the two ordered $9$-tuples of complex numbers $s_X$ and $s_Y$ are coming from the same order-$6$ Somos sequence $s$?

The trouble is that $s_X$ and $s_Y$ can be arbitrarily far apart in the sequence. For any fixed separation between them, we can (in principle) express the elements of one of them as concrete rational functions of the elements of the other, and then compute the desired minor directly. By contrast, the condition that we can run the Somos recurrence some indefinite number of times over $s_X$ so as to obtain $s_Y$ seems difficult to encode algebraically.

So our calculation is not going to rely on the full strength of this condition. Instead, it is only going to make use of some partial information about $s_X$ and $s_Y$. We proceed now to specify what this partial information is. Let $F$ be a rational function of $n$ indeterminates. Then $F$ is a \emph{Somos invariant} of order $n$ when it satisfies \[F(x_0, x_1, \ldots, x_{n - 1}) = F\left(x_1, x_2, \ldots, x_{n - 1}, \frac{a_1x_1x_{n - 1} + a_2x_2x_{n - 2} + \cdots + a_{\lfloor n/2 \rfloor}x_{\lfloor n/2 \rfloor}x_{\lceil n/2 \rceil}}{x_0}\right).\]

Clearly, if $s$ is an order-$n$ Somos sequence with nonzero terms and $F$ is an order-$n$ Somos invariant, then $F$ must remain constant over all windows of size $n$ in $s$. The existence of nontrivial Somos invariants is very much not obvious. Still \cite{H05, H07, H09, FH, U19, TT}, they do exist. We outline one elementary approach to searching for them in Section \ref{inv}.

We say that two Somos sequences $s$ and $t$ are \emph{twins} when they satisfy the same Somos recurrence and all of the relevant Somos invariants (as defined in Section \ref{inv}) agree over them. The desired condition that our calculation hinges on will be simply that $s_X$ and $s_Y$ must be twins. This line of reasoning suggests, furthermore, that Theorem \ref{o6r} ought to admit a generalisation which involves two twinned sequences instead of two copies of the same sequence:

\begin{theorem} \label{o6t} Let $s$ and $t$ be two twinned Somos sequences of order $6$ with nonzero terms. Then, generically, the matrix $s \times t$ is of diamond rank at most $4$. \end{theorem} 

The term ``generically'' in the statement of Theorem \ref{o6t} will be given a precise technical meaning in Section \ref{o6}. Roughly speaking, it tells us that the conclusion holds for ``almost all'' twinned pairs, and the ones for which it does not hold may be viewed as ``degenerate'' in some sense. Notice that we cannot tell yet if Theorem \ref{o6t} implies Theorem \ref{o6r} or not, because we do not know if the twinned pairs of Theorem \ref{o6r} are ``generic'' in the particular way required by Theorem~\ref{o6t}. However, we will see in Section \ref{o6} that Theorem \ref{o6t} is where the heavy lifting takes place, and that the derivation of Theorem \ref{o6r} from it is not too difficult.

We turn to order $7$ now, beginning once again with the decimations:

\begin{proposition} \label{o7d} Let $s$ be a Somos sequence of order $7$ with nonzero terms. Then every decimation of $s$ by an even factor is a Somos sequence of order at most $8$; while every decimation of $s$ by an odd factor is a Somos sequence of order at most $9$. \end{proposition} 

This claim is, as before, a corollary of some finite-rank statement about $s \times s$. However, in this instance the matter is complicated by parity considerations coming into play. It is not true anymore that all $5 \times 5$ diamond minors in $s \times s$ vanish -- though many of them still do. We will introduce, in Section \ref{init}, the notion of a ``half-diamond minor''. Roughly speaking, this is a diamond minor where either the diagonals or the anti-diagonals satisfy a certain parity condition.

\begin{theorem} \label{o7r} Let $s$ be a Somos sequence of order $7$ with nonzero terms. Then the matrix $s \times s$ is of half-diamond rank at most $4$. \end{theorem} 

Both Proposition \ref{o7d} and Theorem \ref{o7r} appear to be new. The same proof method works as with order $6$, and we arrive at an analogous generalisation:

\begin{theorem} \label{o7t} Let $s$ and $t$ be two twinned Somos sequences of order $7$ with nonzero terms. Then, generically, the matrix $s \times t$ is of half-diamond rank at most $4$. \end{theorem} 

We mentioned in the beginning that many of the remarkable properties of Somos sequences do not extend to orders $n \ge 8$. However, for the subclass of Gale-Robinson sequences, some of these properties have been confirmed \cite{G, FZ} to hold universally. We give the definition of this subclass in Section \ref{high}; it includes all Somos sequences of orders $3 \le n \le 7$. (Beyond that, when $n \ge 8$, the subclass and the full class diverge, with the former constituting an ever smaller fraction of the latter.) A conjecture of Ustinov \cite{U24} asserts, in part, that the Gale-Robinson sequences exhibit finite-rank properties, too.

Below, we put forward two conjectures which make more specific predictions. Both of them are certainly true of all orders $3 \le n \le 7$, as shown by Theorems \ref{o6r} and \ref{o7r} in conjunction with the lower-order material of Section \ref{low}. Further evidence in support of Conjectures \ref{gre} and \ref{gro} will be provided in Section \ref{high}.

\begin{conjecture} \label{gre} Let $s$ be a Gale-Robinson sequence with nonzero terms, of a primitive type and even order $n$ with $n = 2m + 2$. Then the matrix $s \times s$ is of diamond rank at most $2^m$. \end{conjecture} 

\begin{conjecture} \label{gro} Let $s$ be a Gale-Robinson sequence with nonzero terms, of a primitive type and odd order $n$ with $n = 2m + 3$. Then the matrix $s \times s$ is of half-diamond rank at most $2^m$. \end{conjecture} 

The rest of the paper is structured as follows: Sections \ref{init} and \ref{algp} cover the basics. Section~\ref{cmin} explains how we make the leap from contiguous minors to arbitrary ones. Section~\ref{inv} introduces the invariants. Sections \ref{o6} and \ref{o7} establish our main results for orders $6$ and $7$, respectively. Section~\ref{int} presents some applications of these results to questions of integrality. Section~\ref{low} reviews the lower-order analogues of the preceding developments. Finally, in Section~\ref{high} we discuss potential higher-order analogues as well.

\section{Initial Observations} \label{init}

To us, a sequence $s$ is a function whose domain is some integer interval $I$. We allow both finite and infinite intervals, and we say that $I$ \emph{indexes} $s$. We often write $s_i$ instead of $s(i)$ for the $i$-th term of $s$.

Let $n$ be a positive integer with $n \ge 2$. Suppose that $s$ is an order-$n$ Somos sequence which satisfies \[s_is_{i + n} = a_1s_{i + 1}s_{i + n - 1} + a_2s_{i + 2}s_{i + n - 2} + \cdots + a_{\lfloor n/2 \rfloor}s_{i + \lfloor n/2 \rfloor}s_{i + \lceil n/2 \rceil} \tag{\textbf{S}}\] for all $i$. (Or, to be precise, for all $i$ such that $[i; i + n] \subseteq I$. From now on, ``for all indices'' will be implicitly understood to mean ``for all valid indices''.) We call $\mathbf{a} = (a_1, a_2, \ldots, a_{\lfloor n/2 \rfloor})$ the \emph{coefficients} of (\textbf{S}), and we denote the class of all order-$n$ Somos sequences which satisfy (\textbf{S}) with the coefficient tuple $\mathbf{a}$ by $\mathfrak{S}_n(\mathbf{a})$.

In this context, we assume by convention that $[0; n - 1] \subseteq I$ and we call $\mathbf{s} = (s_0, s_1, \ldots, s_{n - 1})$ the \emph{seed} of $s$. Clearly, if all terms of $s$ are nonzero, then the coefficients and the seed (together with, implicitly, the indexing interval) determine $s$ uniquely. Conversely, given two ordered tuples of complex numbers $\mathbf{a} \in \mathbb{C}^{\lfloor n/2 \rfloor}$ and $\mathbf{s} \in \mathbb{C}^n$, we can apply the Somos recurrence both forwards and backwards so as to generate an order-$n$ Somos sequence based on them. In both directions, if we run into division by zero, we stop.

The \emph{unit} Somos sequence of order $n$ is the one defined by $\mathbf{a} = (1, 1, \ldots, 1)$ and $\mathbf{s} = (1, 1, \ldots, 1)$, indexed by $\mathbb{Z}$. We never run into division by zero as we construct this sequence because, by induction on the index, all of its terms are positive reals.

Let $c$ be a nonzero complex constant. Notice that, if we multiply each term of $s$ by $c$, the resulting sequence will still be in $\mathfrak{S}_n(\mathbf{a})$. For odd $n$, we can also scale the terms of $s$ parity-wise; i.e., we can multiply all even-indexed terms of $s$ by $c$ while preserving the odd-indexed ones, or vice versa. Furthermore, for both even and odd $n$, multiplying the $i$-th term of $s$ by $c^i$ for all $i$ will once again produce another element of $\mathfrak{S}_n(\mathbf{a})$.

We call these transformations the \emph{symmetries} of order-$n$ Somos sequences. To make sure we are not missing any low-hanging fruit, consider more generally the transformation $s_i \to c^{e_i}s_i$, where $e$ is some as of yet unknown sequence indexed by $\mathbb{Z}$. We stipulate that $e_i + e_{i + n} = e_{i + 1} + e_{i + n - 1} = \cdots = e_{i + \lfloor n/2 \rfloor} + e_{i + \lceil n/2 \rceil}$ for all $i$, so that this transformation is guaranteed to send $\mathfrak{S}_n(\mathbf{a})$ into itself. The complex-valued $e$ which satisfy these constraints form a linear space~$\mathcal{E}$. It is straightforward to see that, if $n$ is even, then $\dim \mathcal{E} = 2$ and one basis for it is given by the sequences $e^\text{I}$ and $e^\text{II}$ defined by $e^\text{I}_i = 1$ and $e^\text{II}_i = i$ for all $i$; whereas, if $n$ is odd, then $\dim \mathcal{E} = 3$ and one basis for it is given by $e^\text{II}$ together with the sequences $e^\text{III}$ and $e^\text{IV}$ defined by $e^\text{III}_i = i \bmod 2$ and $e^\text{IV}_i = (i + 1) \bmod 2$ for all $i$.

We already know that an order-$n$ Somos sequence with nonzero terms is determined uniquely by its coefficients $\mathbf{a}$ and its seed $\mathbf{s}$. We proceed now to make the way in which it is determined by them somewhat more explicit.

Let $\alpha_1$, $\alpha_2$, $\ldots$, $\alpha_{\lfloor n/2 \rfloor}$ and $x_0$, $x_1$, $\ldots$, $x_{n - 1}$ be formal indeterminates. We associate the $\alpha$'s with $a_1$, $a_2$, $\ldots$, $a_{\lfloor n/2 \rfloor}$ and the $x$'s with $s_0$, $s_1$, $\ldots$, $s_{n - 1}$. For convenience, we write $\boldsymbol{\alpha} = (\alpha_1, \alpha_2, \ldots, \alpha_{\lfloor n/2 \rfloor})$ and $\mathbf{x} = (x_0, x_1, \ldots, x_{n - 1})$, similarly to how we defined $\mathbf{a}$ and $\mathbf{s}$ before. Let $\mathcal{A} = \mathbb{Z}[\boldsymbol{\alpha}]$ be the ring of all integer-coefficient polynomials of $\alpha_1$, $\alpha_2$, $\ldots$, $\alpha_{\lfloor n/2 \rfloor}$ and let $\mathcal{A}_\text{Frac} = \mathbb{Z}(\boldsymbol{\alpha})$ be the field of all integer-coefficient rational functions of $\alpha_1$, $\alpha_2$, $\ldots$, $\alpha_{\lfloor n/2 \rfloor}$.

Consider the order-$n$ ``Somos sequence'' $S$ with coefficients $\boldsymbol{\alpha}$ and seed $\mathbf{x}$, indexed by $\mathbb{Z}$. Its terms will be elements of $\mathcal{A}(\mathbf{x})$; i.e., they will be rational functions of $x_0$, $x_1$, $\ldots$, $x_{n - 1}$ with coefficients drawn out of $\mathcal{A}$. We call $S$ the \emph{master} Somos sequence of order $n$. The choice of term is because all complex-number Somos sequences of order $n$ with nonzero terms can be obtained from $S$ by assigning concrete complex numbers to its indeterminates. Explicitly, if all terms of $s$ are nonzero, then $s_i = S_i(\mathbf{a}, \mathbf{s})$ for all $i$. We never run into division by zero as we construct $S$ because we never run into division by zero as we construct the unit Somos sequence.

Let $\sigma$ be the substitution given by \[\sigma(F) = F\left(x_1, x_2, \ldots, x_{n - 1}, \frac{\alpha_1x_1x_{n - 1} + \alpha_2x_2x_{n - 2} + \cdots + \alpha_{\lfloor n/2 \rfloor}x_{\lfloor n/2 \rfloor}x_{\lceil n/2 \rceil}}{x_0}\right)\] for all $F \in \mathcal{A}(\mathbf{x})$. We write $\sigma^i$ for the $i$-th iteration of $\sigma$; since $\sigma$ is invertible, this definition makes sense for negative $i$, too. Then $S_i = \sigma^i(x_0)$ for all $i$.

We also revise our definition of a Somos invariant with the help of $\sigma$. We say that $F$ is a Somos invariant of order $n$ when it is a fixed point of $\sigma$ in $\mathcal{A}(\mathbf{x})$. (Our definition in the introduction was slightly different. In it, the $\alpha$'s of $F$ were substituted with their corresponding $a$'s, for the sake of simplicity.)

We turn now from sequences to matrices. To us, a matrix $M$ is a function whose domain is the product $I \times J$ of two integer intervals $I$ and $J$. We allow both finite and infinite intervals. We say that $I \times J$ \emph{indexes} $M$, and we call its elements the \emph{positions} of $M$. We also call $M(i, j)$ the \emph{entry} of $M$ at position $(i, j)$.

The \emph{diagonal} of $M$ with offset $c$ is the set of all positions $(i, j)$ in $M$ such that $i - j = c$. Similarly, the \emph{anti-diagonal} of $M$ with offset $c$ is defined by $i + j = c$. Each position of $M$ lies on one diagonal and one anti-diagonal. The diagonal with offset $c'$ and the anti-diagonal with offset $c''$ meet at $c' \di c'' = ((c' + c'')/2, (c'' - c')/2)$, provided that the latter is a position of $M$.

We proceed next to formalise the notion of a diamond sub-matrix from the introduction. Let $e'$ and $e''$ be two integer sequences, indexed by $E'$ and $E''$, respectively. For each one of $e'$ and $e''$, we assume that its terms are pairwise distinct. We define the diamond sub-matrix of $M$ with offsets $e'$ and $e''$ to be the matrix indexed by $E' \times E''$ whose entry at position $(i, j)$ equals $M(e'(i) \di e''(j))$. We are only allowed to form this diamond sub-matrix when all of the pairwise intersections $e'(i) \di e''(j)$ are indeed positions of $M$.

For $e'$ and $e''$ to be the offsets of some diamond sub-matrix, the terms of $e'$ and $e''$ must necessarily be all of the same parity. We say that a diamond sub-matrix is \emph{contiguous} when both of $e'$ and $e''$ are arithmetic progressions with common difference $2$. Clearly, this is the ``tightest'' that a diamond sub-matrix can possibly be.

We define a \emph{half-diamond sub-matrix} of $M$ to be a diamond sub-matrix of $M$ where at least one of the offset sequences $e'$ and $e''$ satisfies the stronger condition that all of its terms are congruent modulo $4$. (This definition might seem opaque at first. However, it arises naturally upon careful consideration of the final part of Proposition \ref{o7e} below.) Notice that, in this context, we must revise our definition of contiguity so that it correctly describes the ``tightest'' of these objects. We say that a half-diamond sub-matrix is \emph{contiguous} when one of $e'$ and $e''$ is an arithmetic progression with common difference $2$ and the other one is an arithmetic progression with common difference $4$.

We define the \emph{diamond minor} of $M$ with offsets $e'$ and $e''$ to be the determinant of the diamond sub-matrix of $M$ with offsets $e'$ and $e''$. The notions we just introduced for diamond sub-matrices all carry over in an obvious manner to diamond minors as well. The definition of diamond rank given in the introduction has now been put on firm formal footing. We also define the \emph{half-diamond rank} of $M$ to be the smallest nonnegative integer $r$ such that all half-diamond minors in $M$ of size $(r + 1) \times (r + 1)$ vanish.

The rest of this section sheds some more light on the decimation properties of Somos sequences. This material will not be required for the proofs of Theorems \ref{o6r}--\ref{o7t}.

First we outline how Somos sequences may be approached experimentally. For this, it will be convenient to relax our definition of a Somos sequence somewhat. Suppose that $s$ satisfies the relation \[a_0s_is_{i + n} + a_1s_{i + 1}s_{i + n - 1} + \cdots + a_{\lfloor n/2 \rfloor}s_{i + \lfloor n/2 \rfloor}s_{i + \lceil n/2 \rceil} = 0 \tag{$\mathbf{S}^\star$}\] for all $i$. We waive the condition that $a_0$ must be nonzero, and instead we require merely that the $a$'s are not all zeroes. We say, then, that $s$ is a Somos sequence of nonstrict order $n$. We let $\mathbf{a}^\star = (a_0, a_1, \ldots, a_{\lfloor n/2 \rfloor})$, and we denote the class of all $s$ which satisfy ($\mathbf{S}^\star$) with the coefficient tuple $\mathbf{a}^\star$ by $\mathfrak{S}^\star_n(\mathbf{a}^\star)$.

Of course, if $s$ is of nonstrict order $n$, then it is also of order $m$ for some $m \le n$ of the same parity as $n$. Conversely, if $s$ is of nonstrict order $n$, then certainly it is also of nonstrict order $m$ for all $m \ge n$ of the same parity as $n$. On the other hand, $s$ being of nonstrict order $n$ does not necessarily imply that it is also of nonstrict order $n + 1$. For example, the unit Somos sequence of order $6$ is not a Somos sequence of nonstrict order $7$. So, in particular, ``nonstrict order $n$'' is a distinct concept from ``order at most $n$''.

Suppose that we are given a finite sequence of complex numbers $s$, and we wish to test whether it is a Somos sequence or not. For any fixed nonstrict order $n$, we may go about this task by treating the coefficients $a_0$, $a_1$, $\ldots$, $a_{\lfloor n/2 \rfloor}$ as ``unknowns'', and each instance of ($\mathbf{S}^\star$) as a ``constraint'' imposed upon these unknowns. Taken together, all such constraints form a system of homogeneous linear equations. (The reason we prefer to work with nonstrict orders in this setting is precisely because this makes our linear equations homogeneous.) What remains is to investigate whether this system admits a nontrivial solution. Notice that the matrix of our system will always be a diamond sub-matrix of $s \times s$.

We can now revise Propositions \ref{o6d} and \ref{o7d}, taking the above considerations into account:

\begin{proposition} \label{o6e} Let $s$ be a Somos sequence of order $6$ with nonzero terms. Then every decimation of $s$ is a Somos sequence of nonstrict orders both $8$ and $9$. Furthermore, for each positive integer $d$ and each nonstrict order $n \in \{8, 9\}$, there exist coefficients $\mathbf{a}^\star$ such that all decimations of $s$ by a factor of $d$ belong to the same class $\mathfrak{S}^\star_n(\mathbf{a}^\star)$. \end{proposition} 

\begin{proposition} \label{o7e} Let $s$ be a Somos sequence of order $7$ with nonzero terms. Then every decimation of $s$ by an even factor is a Somos sequence of nonstrict orders both $8$ and $9$, while every decimation of $s$ by an odd factor is a Somos sequence of nonstrict orders both $9$ and $16$. Furthermore, for each positive integer $d$ and each nonstrict order $n$ thus associated with $d$, there exist coefficients $\mathbf{a}^\star$ such that all decimations of $s$ by a factor of $d$ belong to the same class~$\mathfrak{S}^\star_n(\mathbf{a}^\star)$. \end{proposition} 

These nonstrict orders are the best possible, in the sense that none of them can be replaced with smaller ones without the claims becoming false. The unit Somos sequence of order $6$ with decimation factor $2$ confirms this for Proposition \ref{o6e}, while the unit Somos sequence of order~$7$ with decimation factors $2$ and $3$ confirms it for Proposition \ref{o7e}. The revised Propositions \ref{o6e}~and~\ref{o7e} continue to be corollaries of Theorems \ref{o6r} and \ref{o7r}, respectively. The derivations are not too difficult.

\section{Algebraic Preliminaries} \label{algp}

Here, we review some basic notions from algebraic geometry. We develop all of them from scratch, in keeping with our promise of an elementary level of exposition. The purpose of these notions will be to help us deal away with the degeneracies which occur in Sections \ref{o6} and \ref{o7} -- the ``technical fiddling'' referred to in the introduction.

Fix a positive integer $n$. Let $x_0$, $x_1$, $\ldots$, $x_{n - 1}$ be formal indeterminates, and consider the ring $\mathcal{C} = \mathbb{C}[x_0, x_1, \ldots, x_{n - 1}]$ of all complex-coefficient polynomials of $x_0$, $x_1$, $\ldots$, $x_{n - 1}$. We write $\zp$ for the zero polynomial, with $\deg \zp = -\infty$. Below, by a ``system'' we mean any finite set of elements of $\mathcal{C}$, unless otherwise specified.

Let $\mathcal{P}$ be some property which a point $\mathbf{c} = (c_0, c_1, \ldots, c_{n - 1}) \in \mathbb{C}^n$ might or might not possess. We say that $\mathcal{P}$ holds \emph{generically} if there exists a system of nonzero polynomials $\{P_1, P_2, \ldots, P_d\}$ such that $\mathcal{P}$ is true of all $\mathbf{c}$ with $P_1(\mathbf{c}) \neq 0$, $P_2(\mathbf{c}) \neq 0$, $\ldots$, $P_d(\mathbf{c}) \neq 0$.

Intuitively, genericity tells us that $\mathcal{P}$ holds ``almost always'', with the exceptions being ``degenerate'' somehow. The polynomials $P_1$, $P_2$, $\ldots$, $P_d$ serve to describe the various potential degeneracies. Clearly, in the definition of genericity we can assume without loss of generality that the implicit system consists of a single polynomial, by replacing $P_1$, $P_2$, $\ldots$, $P_d$ with their product. Still, it is often more convenient to allow multiple non-degeneracy conditions.

(Notice that we must augment our definition of genericity somehow in the setting of Theorems~\ref{o6t}~and~\ref{o7t}, as the parameters of their statements do not range freely but must instead satisfy certain constraints. These augmentations will be taken care of in Sections \ref{o6} and \ref{o7}.)

For properties $\mathcal{P}$ which can be expressed as ``$P(\mathbf{c}) \neq 0$'', with $P \in \mathcal{C}$, it is typically very easy to show that $\mathcal{P}$ holds generically. Indeed, any concrete $\mathbf{c}$ which satisfies $\mathcal{P}$ would immediately guarantee that $P \neq \zp$. So, when we encounter such properties in Sections \ref{o6} and \ref{o7}, we will usually assert their genericity without further justification.

Sometimes, if we know that a property holds generically, we can conclude on this basis that in fact it holds universally. We say that $\mathcal{P}$ is \emph{algebraic} if there exists a system $\{P_1, P_2, \ldots, P_d\}$ such that $\mathcal{P}$ is true of $\mathbf{c}$ if and only if $P_1(\mathbf{c}) = 0$, $P_2(\mathbf{c}) = 0$, $\ldots$, $P_d(\mathbf{c}) = 0$. Or, in other words, an algebraic property is one which can be expressed as a system of polynomial equations.

\begin{lemma} \label{gen} Suppose that $\mathcal{P}$ is algebraic and that it holds for a generic $\mathbf{c} \in \mathbb{C}^n$. Then it actually holds for all $\mathbf{c} \in \mathbb{C}^n$. \end{lemma} 

\begin{myproof} Let $\{P_1, P_2, \ldots, P_d\}$ be the implicit system we get out of $\mathcal{P}$ being algebraic, and (without loss of generality) let $\{Q\}$ be the implicit system we get out of $\mathcal{P}$ holding generically. Then, for each $i$, the polynomial $P_iQ$ must vanish over all of $\mathbb{C}^n$. So $P_iQ = \zp$. Since $Q \neq \zp$, we conclude that also $P_i = \zp$. \end{myproof}

We go on to an overview of some concrete algebraic properties which will be useful to us later on. Let $y$ be a new formal indeterminate, and consider the ring $\mathcal{C}[y]$ of all polynomials of $y$ whose coefficients are drawn out of $\mathcal{C}$. Let $P$ be any nonzero element of $\mathcal{C}[y]$, with $k = \deg P$. Then, for every $\mathbf{c} \in \mathbb{C}^n$, we get that $P(\mathbf{c}, y)$ is an element of $\mathbb{C}[y]$. Notice that $k$ and $\deg P(\mathbf{c}, y)$ might differ -- even though they coincide generically. When $k = \deg P(\mathbf{c}, y)$, we say that $P(\mathbf{c}, y)$ is of \emph{full degree}.

Let $Q$ be one more nonzero element of $\mathcal{C}[y]$, with $\ell = \deg Q$. The property ``$P(\mathbf{c}, y)$ divides $Q(\mathbf{c}, y)$'' is not always algebraic. For example, if $n = 1$, $P = x_0y$, and $Q = y$, it becomes equivalent to ``$x_0 \neq 0$'', whose non-algebraicity is obvious. However, it is possible to tweak this property very slightly so as to make it algebraic:

\begin{lemma} \label{div} The property ``$P(\mathbf{c}, y)$ divides $Q(\mathbf{c}, y)$, or else $P(\mathbf{c}, y) = \zp$'' is algebraic. \end{lemma} 

\begin{myproof} Let $\lambda_0$, $\lambda_1$, $\ldots$, $\lambda_{\ell + 1}$ be new formal indeterminates. We view the $\lambda$'s as unknowns, and the polynomial equation $(\lambda_0 + \lambda_1y + \cdots + \lambda_\ell y^\ell)P = \lambda_{\ell + 1}Q$ as a constraint imposed upon them. This gives us a system $\Lambda$ of homogeneous linear equations over the $\lambda$'s, whose matrix we denote by $M$. Our desired property of $\mathbf{c}$ is equivalent to $\Lambda(\mathbf{c})$ admitting a nontrivial solution; which, in turn, is equivalent to $M(\mathbf{c})$ being of rank at most $\ell + 1$; which, in turn, is equivalent to $\mathbf{c}$ being a root of every minor in $M$ of size $(\ell + 2) \times (\ell + 2)$. \end{myproof}

Consider, next, the property ``$P(\mathbf{c}, y)$ and $Q(\mathbf{c}, y)$ share a non-constant common factor''. It is not quite algebraic, either. For example, if $n = 2$, $P = x_0y + 1$, and $Q = x_1y + 1$, it becomes equivalent to ``$x_0 = x_1 \neq 0$'', whose non-algebraicity is straightforward. Once again, though, we can patch things up without too much trouble:

\begin{lemma} \label{cf} The property ``$P(\mathbf{c}, y)$ and $Q(\mathbf{c}, y)$ share a non-constant common factor, or else neither one of them is of full degree'' is algebraic. \end{lemma} 

\begin{myproof} The argument is similar to the one we employed in the proof of Lemma~\ref{div}. The case when $k = \ell = 0$ is trivial. Otherwise, let $\lambda_0$, $\lambda_1$, $\ldots$, $\lambda_{k + \ell - 1}$ be new formal indeterminates. This time around, we extract our system of homogeneous linear equations $\Lambda$ out of the polynomial equation $(\lambda_0 + \lambda_1y + \cdots + \lambda_{\ell - 1}y^{\ell - 1})P = (\lambda_\ell + \lambda_{\ell + 1}y + \cdots + \lambda_{k + \ell - 1}y^{k - 1})Q$. It is associated with a square matrix $M$ of size $(k + \ell) \times (k + \ell)$. So, in the present setting, our desired property becomes equivalent to $\mathbf{c}$ being a root of $\det M$. \end{myproof}

The polynomial $\det M$ which arises in the proof is known as the \emph{resultant} of $P$ and $Q$ with respect to $y$. We denote it by $\res_y(P, Q)$. (This definition assumes that $k + \ell \ge 1$.) Notice also that the proof of Lemma \ref{cf} continues to hold when $\mathbb{C}$ is replaced with an arbitrary field. This observation will be important in Sections \ref{o6} and \ref{o7}.

Recall that an element of $\mathbb{C}[y]$ is divisible by a non-constant square if and only if it shares a non-constant common factor with its formal derivative. The resultant $\res_y(P, \partial_y P)$ is known as the \emph{discriminant} of $P$ with respect to $y$, and we denote it by $\disc_y(P)$. (This definition assumes that $k \ge 1$.) We arrive at the following corollary of Lemma \ref{cf}:

\begin{lemma} \label{sf} The property ``$P(\mathbf{c}, y)$ is divisible by a non-constant square, or else it is not of full degree'' is algebraic. \end{lemma} 

We will not refer to Lemmas \ref{cf} and \ref{sf} explicitly in future sections. Instead, we will refer to the relevant resultants and discriminants directly.

\section{Contiguous Minors} \label{cmin}

Let $M$ be a matrix over any field.

\begin{lemma} \label{cm} Suppose that all contiguous minors in $M$ of size $(r + 1) \times (r + 1)$ vanish. Suppose also that no contiguous minors in $M$ of size $r \times r$ do. Then $M$ is of rank at most $r$. So, in particular, all minors in $M$ of size $(r + 1) \times (r + 1)$ vanish, the non-contiguous ones included. \end{lemma} 

Of course, if $M$ is of size at least $r \times r$, so that the non-vanishing condition is not vacuous, we can omit the ``at most''.

\begin{myproof} The claim is trivial when $M$ contains at most $r$ rows or at most $r$ columns. Suppose, from now on, that $M$ is of size at least $(r + 1) \times (r + 1)$.

We begin with the special case when $M$ consists of $r + 1$ rows exactly. For all $i \le j$, let $K(i, j)$ be the linear hull of columns $i$, $i + 1$, $\ldots$, $j$ of $M$. Since columns $i$, $i + 1$, $\ldots$, $i + r$ form a contiguous sub-matrix in $M$ of size $(r + 1) \times (r + 1)$, we get that $\dim K(i, i + r) \le r$. On the other hand, the union of columns $i + 1$, $i + 2$, $\ldots$, $i + r$ contains a contiguous sub-matrix in $M$ of size $r \times r$, and so $\dim K(i + 1, i + r) = r$. But clearly $K(i + 1, i + r) \subseteq K(i, i + r)$, and so $K(i, i + r) = K(i + 1, i + r)$. By the same token, also $K(i + 1, i + r) = K(i + 1, i + r + 1)$. We get that $K(i, i + r) = K(i + 1, i + r + 1)$ for all $i$. So, as $i$ varies, $K(i, i + r)$ remains constant. Denote its constant value by $K$. Then all columns of $M$ are in $K$ and, as previously noted, $\dim K = r$.

For the general case, we run the same argument one more time, but ``vertically''. Suppose that $M$ contains at least $r + 2$ rows. For all $i \le j$, let $L(i, j)$ be the linear hull of rows $i$, $i + 1$, $\ldots$, $j$ of $M$. Since rows $i$, $i + 1$, $\ldots$, $i + r$ form a sub-matrix in $M$ of height $r + 1$ which satisfies the conditions of Lemma \ref{cm}, by the preceding discussion we find that $\dim L(i, i + r) = r$. On the other hand, the union of rows $i + 1$, $i + 2$, $\ldots$, $i + r$ contains a contiguous sub-matrix in $M$ of size $r \times r$, and so $\dim L(i + 1, i + r) = r$ as well. But clearly $L(i + 1, i + r) \subseteq L(i, i + r)$, and so once again we arrive at $L(i, i + r) = L(i + 1, i + r)$. By the same token, also $L(i + 1, i + r) = L(i + 1, i + r + 1)$. The rest of the argument goes as before. \end{myproof}

Notice that we cannot afford to miss even a single contiguous minor of size $(r + 1) \times (r + 1)$. For example, consider the matrix $M_\#$ over $\mathbb{Q}$, indexed by $\mathbb{Z} \times \mathbb{Z}$, which we form as follows: First, fill all positions with zeroes. Then, at all positions $(i, j)$ with $i \equiv j \pmod r$, replace the $0$ with~a~$1$. Finally, at all positions $(i, j)$ with $i \ge 0$, $j \ge 0$, and $i \equiv j \equiv 0 \pmod r$, replace the $1$ with~a~$2$. It is straightforward to see that in $M_\#$ all contiguous minors of size $r \times r$ are nonzero, while all but one contiguous minors of size $(r + 1) \times (r + 1)$ vanish.

So the vanishing condition cannot be weakened. The non-vanishing condition, though, is a different matter. It turns out that we do not need to inspect every single contiguous minor of size $r \times r$ in $M$ so as to confirm that they are all non-vanishing. We can instead get away with inspecting just a small fraction of them. Since such optimisations are not crucial to our main task (of proving Theorems \ref{o6r}--\ref{o7t}), we limit ourselves to some simple observations.

\begin{proposition} \label{dia} Let $M$ be a matrix indexed by $\mathbb{Z} \times \mathbb{Z}$. Suppose that all contiguous minors in $M$ of size $(r + 1) \times (r + 1)$ vanish. Consider the contiguous minors in $M$ of size $r \times r$ ``on the main diagonal''; i.e., the ones where the integer intervals which index the rows and the columns coincide. Suppose also that none of them vanish. Then no contiguous minors in $M$ of size $r \times r$ vanish at all. So, in particular, $M$ satisfies the conditions of Lemma \ref{cm} as well as its conclusion. \end{proposition} 

\begin{myproof} Let $M_\P$ be the matrix whose entry at position $(i, j)$ equals the contiguous $r \times r$ minor in $M$ defined by columns $i$, $i + 1$, $\ldots$, $i + r - 1$ and rows $j$, $j + 1$, $\ldots$, $j + r - 1$.

Consider any matrix $W$ of size $(r + 1) \times (r + 1)$. There are four contiguous sub-matrices of size $r \times r$ in $W$. Denote the ones in top left, top right, lower left, lower right by $W_1$, $W_2$, $W_3$, $W_4$, respectively. Let also $W_\S$ be the unique contiguous sub-matrix of size ${(r - 1)} \times {(r - 1)}$ in $W$ which is concentric with $W$. The well-known Desnanot–Jacobi identity states that $\det W_1 \det W_4 - \det W_2 \det W_3 = \det W_\S \det W$.

In the setting of Proposition \ref{dia}, this tells us that every contiguous sub-matrix $\big(\begin{smallmatrix} w_1 & w_2\\ w_3 & w_4 \end{smallmatrix}\big)$ of size $2 \times 2$ in $M_\P$ satisfies $w_1w_4 = w_2w_3$. Thus, if any diagonal in $M_\P$ consists entirely of nonzero entries, then so must both of its neighbouring diagonals as well. By induction on the offset, with the main diagonal of $M_\P$ as our base case, we conclude that in fact all entries of $M_\P$ must be nonzero, as desired. \end{myproof}

The same argument shows that, in a finite $M$ of size $(m + r - 1) \times (n + r - 1)$, it suffices to inspect just $\max\{m, n\}$ contiguous minors of size $r \times r$ instead of all $mn$ of them. Other configurations work as well. For example, in any $M$, finite or infinite, consider any cross $C$ formed as the union of one strip of $r$ successive rows and one strip of $r$ successive columns. By the same reasoning as in the proof of Proposition \ref{dia}, we get that it suffices to inspect just the contiguous minors of size $r \times r$ contained within $C$.

\section{The Invariants} \label{inv}

Fix a positive integer $n \ge 2$, and let $\Pi = x_0x_1 \cdots x_{n - 1}$. We will be looking for order-$n$ Somos invariants of the particular form $F = \Phi/\Pi$, with $\Phi$ being a homogeneous polynomial of degree $n$ in $\mathcal{A}[\mathbf{x}]$.

The homogeneous polynomials of degree $n$ in $\mathcal{A}_\text{Frac}[\mathbf{x}]$ form a linear space $\Upsilon$ over the field $\mathcal{A}_\text{Frac}$ in which the unit-coefficient monomials constitute a basis. Consider the transformation \begin{multline*} \varphi(\Phi) = x_0^{n - 2}(\alpha_1x_1x_{n - 1} + \alpha_2x_2x_{n - 2} + \cdots + \alpha_{\lfloor n/2 \rfloor}x_{\lfloor n/2 \rfloor}x_{\lceil n/2 \rceil})\Phi(x_0, x_1, \ldots, x_{n - 1}) - {}\\ {} - \Phi(x_0x_1, x_0x_2, \ldots, x_0x_{n - 1}, \alpha_1x_1x_{n - 1} + \alpha_2x_2x_{n - 2} + \cdots + \alpha_{\lfloor n/2 \rfloor}x_{\lfloor n/2 \rfloor}x_{\lceil n/2 \rceil}) \end{multline*} over $\Upsilon$. It is straightforward to see that $\varphi$ is linear, and that $F$ is a Somos invariant if and only if $\Phi$ belongs to the kernel of $\varphi$.

So, in order to find all Somos invariants of our desired form, it suffices to compute this kernel. Before we get around to that, though, we are going to impose one additional constraint on $F$~and~$\Phi$.

Recall the linear space $\mathcal{E}$ of Section \ref{init} which captures the symmetries of order-$n$ Somos sequences. We require that $F$ agrees with $\mathcal{E}$, in the sense that $F(x_0, x_1, \ldots, x_{n - 1}) = F(y^{e_0}x_0, y^{e_1}x_1,\allowbreak \ldots, y^{e_{n - 1}}x_{n - 1})$ for all integer $e \in \mathcal{E}$, with $y$ being a new formal indeterminate. This is equivalent to each exponent tuple $(d_0, d_1, \ldots, d_{n - 1})$ which occurs in $\Phi$ satisfying $d_0e_0 + d_1e_1 + \cdots + d_{n - 1}e_{n - 1} = e_0 + e_1 + \cdots + e_{n - 1}$ for all integer $e \in \mathcal{E}$.

The polynomials $\Phi$ which satisfy our additional constraint form a linear subspace $\Upsilon_\boxtimes$ of $\Upsilon$. So, from now on, we may focus on computing the kernel of $\varphi$ solely over $\Upsilon_\boxtimes$. The additional constraint serves a twofold purpose. First, it makes the computations a lot more manageable. For example, with $n = 6$, it brings the dimension of our linear space from $\dim \Upsilon = 462$ down to $\dim \Upsilon_\boxtimes = 32$; or, with $n = 7$, from $\dim \Upsilon = 1716$ down to $\dim \Upsilon_\boxtimes = 40$. Second, it helps us pin down just the invariants which will be relevant to our purposes. Indeed, for both orders $5$~and~$7$, there exist additional Somos invariants of our desired form where $\Phi$ is in $\Upsilon$ but not in $\Upsilon_\boxtimes$; however, we do not require these invariants for the proofs of Theorems \ref{o7r} and \ref{o7t} or their order-$5$ analogues Theorems \ref{o5r} and \ref{o5t}.

Let $\Omega_\boxtimes$ be the kernel of $\varphi$ over $\Upsilon_\boxtimes$. Notice that $\Pi$ is always in $\Omega_\boxtimes$. It corresponds to the trivial Somos invariant where $F$ is the constant unity. Notice also that our computations use coefficients in $\mathcal{A}_\text{Frac}$, whereas ultimately we want the coefficients of $F$ and $\Phi$ to be in $\mathcal{A}$. We resolve this issue simply by clearing the denominators.

We proceed now to report the results of the computations. We cover all orders $2 \le n \le 7$. (The lower-order invariants will play a key role in Section \ref{low}.) Since we will be referring to many different orders $n$, below we rename $\Pi$ to $\Pi_n$.

For orders $2$ and $3$, we get that $\dim \Omega_\boxtimes = 1$. So we do not obtain any nontrivial invariants.

For order $4$, we get that $\dim \Upsilon_\boxtimes = 5$ and $\dim \Omega_\boxtimes = 2$. So we obtain, in essence, a single nontrivial invariant. We set $\Phi_4$ to \[x_0^2x_3^2 + \alpha_1x_0x_2^3 + \alpha_1x_1^3x_3 + \alpha_2x_1^2x_2^2,\] and we denote this invariant by $F_4 = \Phi_4/\Pi_4$.

For order $5$, we get that $\dim \Upsilon_\boxtimes = 6$ and $\dim \Omega_\boxtimes = 2$ once again. So, just as with order $4$, we obtain a single nontrivial invariant. We set $\Phi_5$ to \[x_0^2x_3^2x_4 + x_0x_1^2x_4^2 + \alpha_1x_0x_2^2x_3^2 + \alpha_1x_1^2x_2^2x_4 + \alpha_2x_1x_2^3x_3,\] and we denote this invariant by $F_5 = \Phi_5/\Pi_5$.

The higher orders yield polynomials with a large number of summands. For the sake of clarity, we will present these polynomials also in tabular form. For each summand in one of them, one row of the corresponding table will list its exponent tuple and its coefficient. To save space, the exponent tuples will be encoded as decimal strings. For example, the summand $\alpha_3x_0^2x_2x_4^2x_5$ of $\Phi_6$ will be represented by the decimal string $201021$, encoding its exponent tuple $(2, 0, 1, 0, 2, 1)$, together with its coefficient $\alpha_3$.

\begin{table} \footnotesize \null \hfill \parbox{.375\textwidth}{\centering \begin{tabular}{cc} $012210$ & $\alpha_3^2$\\ $013020$ & $\alpha_2\alpha_3$\\ $020310$ & $\alpha_2\alpha_3$\\ $021120$ & $\alpha_2^2$\\ $022011$ & $\alpha_1\alpha_3$\\ $030111$ & $\alpha_1\alpha_2$\\ $102201$ & $\alpha_1\alpha_3$\\ $103011$ & $\alpha_1\alpha_2$\\ $110220$ & $\alpha_1\alpha_3$\\ $110301$ & $\alpha_1\alpha_2$\\ $111030$ & $\alpha_1\alpha_2$\\ $120102$ & $\alpha_3$\\ $201021$ & $\alpha_3$\\ $201102$ & $\alpha_2$ \end{tabular} \caption{} \label{i6a}} \hfill \parbox{.375\textwidth}{\centering \begin{tabular}{cc} $003300$ & $\alpha_1\alpha_3^2$\\ $004110$ & $\alpha_1\alpha_2\alpha_3$\\ $011400$ & $\alpha_1\alpha_2\alpha_3$\\ $012210$ & $\alpha_1\alpha_2^2$\\ $013101$ & $\alpha_1^2\alpha_3$\\ $021201$ & $\alpha_1^2\alpha_2 + \alpha_3^2$\\ $022011$ & $\alpha_2\alpha_3$\\ $031002$ & $\alpha_1\alpha_3$\\ $101310$ & $\alpha_1^2\alpha_3$\\ $102120$ & $\alpha_1^2\alpha_2 + \alpha_3^2$\\ $102201$ & $\alpha_2\alpha_3$\\ $110220$ & $\alpha_2\alpha_3$\\ $112002$ & $\alpha_1\alpha_2$\\ $120021$ & $\alpha_1\alpha_2$\\ $200130$ & $\alpha_1\alpha_3$\\ $200211$ & $\alpha_1\alpha_2$\\ $210012$ & $\alpha_3$ \end{tabular} \caption{} \label{i6b}} \hfill \null \end{table}

For order $6$, we get that $\dim \Upsilon_\boxtimes = 32$ and $\dim \Omega_\boxtimes = 3$. So, in this case, we obtain two linearly independent nontrivial invariants. We set $\Phi_6$ and $\Psi_6$ to {\allowdisplaybreaks \begin{gather*} \alpha_2x_0^2x_2x_3x_5^2 + \alpha_3x_0^2x_2x_4^2x_5 + \alpha_3x_0x_1^2x_3x_5^2 + \alpha_1\alpha_2x_0x_1x_2x_4^3 + \alpha_1\alpha_2x_0x_1x_3^3x_5 + {}\\* \alpha_1\alpha_3x_0x_1x_3^2x_4^2 + \alpha_1\alpha_2x_0x_2^3x_4x_5 + \alpha_1\alpha_3x_0x_2^2x_3^2x_5 + \alpha_1\alpha_2x_1^3x_3x_4x_5 + \alpha_1\alpha_3x_1^2x_2^2x_4x_5 + {}\\* \alpha_2^2x_1^2x_2x_3x_4^2 + \alpha_2\alpha_3x_1^2x_3^3x_4 + \alpha_2\alpha_3x_1x_2^3x_4^2 + \alpha_3^2x_1x_2^2x_3^2x_4 \end{gather*} \aftergroup \ignorespaces} and {\allowdisplaybreaks \begin{gather*} \alpha_3x_0^2x_1x_4x_5^2 + \alpha_1\alpha_2x_0^2x_3^2x_4x_5 + \alpha_1\alpha_3x_0^2x_3x_4^3 + \alpha_1\alpha_2x_0x_1^2x_4^2x_5 + \alpha_1\alpha_2x_0x_1x_2^2x_5^2 + {}\\* \alpha_2\alpha_3x_0x_1x_3^2x_4^2 + \alpha_2\alpha_3x_0x_2^2x_3^2x_5 + (\alpha_1^2\alpha_2 + \alpha_3^2)x_0x_2^2x_3x_4^2 + \alpha_1^2\alpha_3x_0x_2x_3^3x_4 + \alpha_1\alpha_3x_1^3x_2x_5^2 + {}\\ \alpha_2\alpha_3x_1^2x_2^2x_4x_5 + (\alpha_1^2\alpha_2 + \alpha_3^2)x_1^2x_2x_3^2x_5 + \alpha_1^2\alpha_3x_1x_2^3x_3x_5 + \alpha_1\alpha_2^2x_1x_2^2x_3^2x_4 + \alpha_1\alpha_2\alpha_3x_1x_2x_3^4 + {}\\* \alpha_1\alpha_2\alpha_3x_2^4x_3x_4 + \alpha_1\alpha_3^2x_2^3x_3^3, \end{gather*} \aftergroup \ignorespaces} respectively; the same polynomials are shown also in Tables \ref{i6a} and \ref{i6b}. We denote the invariants associated with them by $F_6 = \Phi_6/\Pi_6$ and $G_6 = \Psi_6/\Pi_6$.

\begin{table} \footnotesize \null \hfill \parbox{.375\textwidth}{\centering \begin{tabular}{cc} $0023200$ & $\alpha_1\alpha_3^2$\\ $0032110$ & $\alpha_1\alpha_2\alpha_3$\\ $0112300$ & $\alpha_1\alpha_2\alpha_3$\\ $0121210$ & $\alpha_1\alpha_2^2$\\ $0122101$ & $\alpha_1^2\alpha_3 + \alpha_3^2$\\ $0131011$ & $\alpha_2\alpha_3$\\ $0211201$ & $\alpha_1^2\alpha_2$\\ $0221002$ & $\alpha_1\alpha_3$\\ $1012210$ & $\alpha_1^2\alpha_3 + \alpha_3^2$\\ $1013101$ & $\alpha_2\alpha_3$\\ $1021120$ & $\alpha_1^2\alpha_2$\\ $1101310$ & $\alpha_2\alpha_3$\\ $1112002$ & $\alpha_1\alpha_2$\\ $1120021$ & $\alpha_1\alpha_2$\\ $1200211$ & $\alpha_1\alpha_2$\\ $1210012$ & $\alpha_3$\\ $2001220$ & $\alpha_1\alpha_3$\\ $2002111$ & $\alpha_1\alpha_2$\\ $2100121$ & $\alpha_3$ \end{tabular} \caption{} \label{i7a}} \hfill \parbox{.375\textwidth}{\centering \begin{tabular}{cc} $0032110$ & $\alpha_1\alpha_3^2$\\ $0041020$ & $\alpha_1\alpha_2\alpha_3$\\ $0112300$ & $\alpha_1\alpha_3^2$\\ $0121210$ & $2\alpha_1\alpha_2\alpha_3$\\ $0130120$ & $\alpha_1\alpha_2^2$\\ $0131011$ & $\alpha_1^2\alpha_3$\\ $0201400$ & $\alpha_1\alpha_2\alpha_3$\\ $0210310$ & $\alpha_1\alpha_2^2$\\ $0211201$ & $\alpha_1^2\alpha_3 + \alpha_3^2$\\ $0220111$ & $\alpha_1^2\alpha_2 + \alpha_2\alpha_3$\\ $0300301$ & $\alpha_1^2\alpha_2$\\ $0310102$ & $\alpha_1\alpha_3$\\ $1013101$ & $\alpha_1^2\alpha_3$\\ $1021120$ & $\alpha_1^2\alpha_3 + \alpha_3^2$\\ $1022011$ & $\alpha_1^2\alpha_2 + \alpha_2\alpha_3$\\ $1030030$ & $\alpha_1^2\alpha_2$\\ $1101310$ & $\alpha_1^2\alpha_3$\\ $1102201$ & $\alpha_1^2\alpha_2 + \alpha_2\alpha_3$\\ $1110220$ & $\alpha_1^2\alpha_2 + \alpha_2\alpha_3$\\ $1112002$ & $\alpha_1\alpha_3$\\ $1120021$ & $\alpha_1\alpha_3$\\ $1200211$ & $\alpha_1\alpha_3$\\ $1201102$ & $\alpha_1\alpha_2$\\ $2002111$ & $\alpha_1\alpha_3$\\ $2003002$ & $\alpha_1\alpha_2$\\ $2010130$ & $\alpha_1\alpha_3$\\ $2011021$ & $\alpha_1\alpha_2$\\ $2101012$ & $\alpha_3$ \end{tabular} \caption{} \label{i7b}} \hfill \null \end{table}

For order $7$, we get that $\dim \Upsilon_\boxtimes = 40$ and $\dim \Omega_\boxtimes = 3$ once again. So, just as with order $6$, we obtain two linearly independent nontrivial invariants. We set $\Phi_7$ and $\Psi_7$ to {\allowdisplaybreaks \begin{gather*} \alpha_3x_0^2x_1x_4x_5^2x_6 + \alpha_1\alpha_2x_0^2x_3^2x_4x_5x_6 + \alpha_1\alpha_3x_0^2x_3x_4^2x_5^2 + \alpha_3x_0x_1^2x_2x_5x_6^2 + \alpha_1\alpha_2x_0x_1^2x_4^2x_5x_6 + {}\\* \alpha_1\alpha_2x_0x_1x_2^2x_5^2x_6 + \alpha_1\alpha_2x_0x_1x_2x_3^2x_6^2 + \alpha_2\alpha_3x_0x_1x_3x_4^3x_5 + \alpha_1^2\alpha_2x_0x_2^2x_3x_4x_5^2 + {}\\ \alpha_2\alpha_3x_0x_2x_3^3x_4x_6 + (\alpha_1^2\alpha_3 + \alpha_3^2)x_0x_2x_3^2x_4^2x_5 + \alpha_1\alpha_3x_1^2x_2^2x_3x_6^2 + \alpha_1^2\alpha_2x_1^2x_2x_3x_4^2x_6 + {}\\ \alpha_2\alpha_3x_1x_2^3x_3x_5x_6 + (\alpha_1^2\alpha_3 + \alpha_3^2)x_1x_2^2x_3^2x_4x_6 + \alpha_1\alpha_2^2x_1x_2^2x_3x_4^2x_5 + \alpha_1\alpha_2\alpha_3x_1x_2x_3^2x_4^3 + {}\\* \alpha_1\alpha_2\alpha_3x_2^3x_3^2x_4x_5 + \alpha_1\alpha_3^2x_2^2x_3^3x_4^2 \end{gather*} \aftergroup \ignorespaces} and {\allowdisplaybreaks \begin{gather*} \alpha_3x_0^2x_1x_3x_5x_6^2 + \alpha_1\alpha_2x_0^2x_2x_3x_5^2x_6 + \alpha_1\alpha_3x_0^2x_2x_4x_5^3 + \alpha_1\alpha_2x_0^2x_3^3x_6^2 + \alpha_1\alpha_3x_0^2x_3^2x_4x_5x_6 + {}\\* \alpha_1\alpha_2x_0x_1^2x_3x_4x_6^2 + \alpha_1\alpha_3x_0x_1^2x_4^2x_5x_6 + \alpha_1\alpha_3x_0x_1x_2^2x_5^2x_6 + \alpha_1\alpha_3x_0x_1x_2x_3^2x_6^2 + {}\\ (\alpha_1^2\alpha_2 + \alpha_2\alpha_3)x_0x_1x_2x_4^2x_5^2 + (\alpha_1^2\alpha_2 + \alpha_2\alpha_3)x_0x_1x_3^2x_4^2x_6 + \alpha_1^2\alpha_3x_0x_1x_3x_4^3x_5 + \alpha_1^2\alpha_2x_0x_2^3x_5^3 + {}\\ (\alpha_1^2\alpha_2 + \alpha_2\alpha_3)x_0x_2^2x_3^2x_5x_6 + (\alpha_1^2\alpha_3 + \alpha_3^2)x_0x_2^2x_3x_4x_5^2 + \alpha_1^2\alpha_3x_0x_2x_3^3x_4x_6 + \alpha_1\alpha_3x_1^3x_2x_4x_6^2 + {}\\ \alpha_1^2\alpha_2x_1^3x_4^3x_6 + (\alpha_1^2\alpha_2 + \alpha_2\alpha_3)x_1^2x_2^2x_4x_5x_6 + (\alpha_1^2\alpha_3 + \alpha_3^2)x_1^2x_2x_3x_4^2x_6 + \alpha_1\alpha_2^2x_1^2x_2x_4^3x_5 + {}\\ \alpha_1\alpha_2\alpha_3x_1^2x_3x_4^4 + \alpha_1^2\alpha_3x_1x_2^3x_3x_5x_6 + \alpha_1\alpha_2^2x_1x_2^3x_4x_5^2 + {}\\* 2\alpha_1\alpha_2\alpha_3x_1x_2^2x_3x_4^2x_5 + \alpha_1\alpha_3^2x_1x_2x_3^2x_4^3 + \alpha_1\alpha_2\alpha_3x_2^4x_3x_5^2 + \alpha_1\alpha_3^2x_2^3x_3^2x_4x_5, \end{gather*} \aftergroup \ignorespaces} respectively; the same polynomials are shown also in Tables \ref{i7a} and \ref{i7b}. We denote the invariants associated with them by $F_7 = \Phi_7/\Pi_7$ and $G_7 = \Psi_7/\Pi_7$.

\section{Order 6} \label{o6}

Here, we prove Theorems \ref{o6r} and \ref{o6t}. Let $s$ and $t$ be two order-$6$ Somos sequences, both with coefficients $\mathbf{a}$, and with seeds $\mathbf{s}$ and $\mathbf{t}$, respectively. For this section, we specialise the general notations and definitions of Section \ref{init} to order $6$; for example, $S$ is going to denote the master Somos sequence of order $6$ throughout.

Since we will be working with two sequences simultaneously, we must ``clone'' $S$ as well as the order-$6$ invariants obtained in Section \ref{inv}. Let $y_0$, $y_1$, $\ldots$, $y_5$ be new formal indeterminates and let $T$ be the sequence obtained from $S$ by substituting all $x$'s with their corresponding $y$'s. For this section, we rename $\Phi_6$, $\Psi_6$, $\Pi_6$ to $\Phi_X$, $\Psi_X$, $\Pi_X$; and we obtain $\Phi_Y$, $\Psi_Y$, $\Pi_Y$ from them by substituting all $x$'s with their corresponding $y$'s once again. The definitions of $F_X$, $F_Y$, $G_X$, $G_Y$ are analogous.

Let $U$ be the numerator of $F_X - F_Y$; explicitly, $U = \Pi_Y\Phi_X - \Pi_X\Phi_Y$. Similarly, let $V$ be the numerator of $G_X - G_Y$; explicitly, $V = \Pi_Y\Psi_X - \Pi_X\Psi_Y$. (The denominators of both differences equal $\Pi_X\Pi_Y$.) Then $s$ and $t$ are twins if and only if $U(\mathbf{a}, \mathbf{s}, \mathbf{t}) = 0$ and $V(\mathbf{a}, \mathbf{s}, \mathbf{t}) = 0$. We call an $(\mathbf{a}, \mathbf{s}, \mathbf{t})$ with this property \emph{admissible}, and we say that $s$ and $t$ are \emph{based} on $(\mathbf{a}, \mathbf{s}, \mathbf{t})$. For convenience, we view $(\mathbf{a}, \mathbf{s}, \mathbf{t})$ as an element of $\mathbb{C}^{15}$ rather than as an element of $\mathbb{C}^3 \times \mathbb{C}^6 \times \mathbb{C}^6$.

For Theorem \ref{o6t}, we wish to show that all diamond minors of size $5 \times 5$ in $s \times t$ vanish. We begin with the contiguous ones among them, as advertised in the introduction.

\begin{lemma} \label{o6cm} Suppose that $s$ and $t$ are twins with nonzero terms for which $a_2$ and $a_3$ are not both zero. Then all contiguous diamond minors of size $5 \times 5$ in the matrix $s \times t$ vanish. \end{lemma} 

The non-degeneracy condition $(a_2, a_3) \neq (0, 0)$ cannot be omitted. It is in some sense the weakest condition on $\mathbf{a}$ which makes the statement true, as evidenced by the fact that setting $(\alpha_2, \alpha_3) = (0, 0)$ in the argument below causes $U$ and $V$ to vanish but not $D$.

\begin{myproof} It suffices to consider the case when both of $s$ and $t$ are of size $9$, indexed by $[-1; 7]$. Then $s \times t$ contains a unique contiguous diamond minor of size $5 \times 5$, which we denote by $\delta$. Consider next the contiguous diamond minor $\Delta$ of size $5 \times 5$ in $S \times T$ with offsets $(-4, -2, 0, 2, 4)$ and $(2, 4, 6, 8, 10)$. Of course, $\Delta \in \mathcal{A}(\mathbf{x}, \mathbf{y})$ and $\delta = \Delta(\mathbf{a}, \mathbf{s}, \mathbf{t})$.

Let $\Delta = D/x_0y_0\Pi_X\Pi_Y$, with $D \in \mathcal{A}[\mathbf{x}, \mathbf{y}]$. (The denominator of $\Delta$ in lowest terms is in fact $x_0^2x_1x_5y_0^2y_1y_5$. The point of rewriting $\Delta$ in this way is so that its denominator becomes ``synchronised'' with the denominators of $F_X - F_Y$ and $G_X - G_Y$. The numerator $D$ is homogeneous of degree $24$, with $687$ summands.)

We aim to show that $U(\mathbf{a}, \mathbf{s}, \mathbf{t}) = 0$ and $V(\mathbf{a}, \mathbf{s}, \mathbf{t}) = 0$ together force $D(\mathbf{a}, \mathbf{s}, \mathbf{t}) = 0$, subject to the condition that $(a_2, a_3) \neq (0, 0)$. We claim that there exist $A_2$ and $B_2$ in $\mathcal{A}[\mathbf{x}, \mathbf{y}]$ with $A_2U + B_2V = \alpha_2D$; as well as $A_3$ and $B_3$ in $\mathcal{A}[\mathbf{x}, \mathbf{y}]$ with $A_3U + B_3V = \alpha_3D$.

\begin{table}[ht] \footnotesize \null \hfill \parbox{.375\textwidth}{\centering \begin{tabular}{cc} $003300102120$ & $\alpha_2\alpha_3$\\ $003300111111$ & $\alpha_1\alpha_3$\\ $004110102120$ & $\alpha_2^2$\\ $004110111111$ & $\alpha_1\alpha_2$\\ $013101102120$ & $\alpha_1\alpha_2$\\ $013101111111$ & $\alpha_1^2$\\ $021201111111$ & $\alpha_3$\\ $022011111111$ & $\alpha_2$\\ $031002111111$ & $\alpha_1$\\ $101310011400$ & $\alpha_2\alpha_3$\\ $101310012210$ & $\alpha_2^2$\\ $101310021201$ & $\alpha_1\alpha_2$\\ $101310102120$ & $\alpha_1\alpha_2$\\ $101310111111$ & $\alpha_1^2$\\ $101310200211$ & $\alpha_2$\\ $102120111111$ & $\alpha_3$\\ $200130102120$ & $\alpha_2$\\ $200130111111$ & $\alpha_1$\\ $210012111111$ & $1$ \end{tabular} \caption{} \label{o6a2}} \hfill \parbox{.375\textwidth}{\centering \begin{tabular}{cc} $102120102201$ & $\alpha_2$\\ $102120110220$ & $\alpha_2$\\ $110301101310$ & $\alpha_2$\\ $111111012210$ & $\alpha_3$\\ $111111013020$ & $\alpha_2$\\ $111111022011$ & $\alpha_1$\\ $111111102201$ & $\alpha_1$\\ $111111110220$ & $\alpha_1$\\ $111111120102$ & $1$\\ $111111201021$ & $1$ \end{tabular} \caption{} \label{o6b2}} \hfill \null \end{table}

The $A$'s and $B$'s we are about to present exhibit certain symmetries, and we will exploit these symmetries so as to compress the presentations. We define the \emph{skew-symmetrisation} of any $P \in \mathcal{A}[\mathbf{x}, \mathbf{y}]$ to be $P(\mathbf{x}, \mathbf{y}) - P(\mathbf{y}, \mathbf{x})$. Our $A_2$ and $B_2$ are the skew-symmetrisations of the polynomials shown in Tables \ref{o6a2} and \ref{o6b2}, respectively. (The exponents of the $x$'s and the $y$'s are listed in the order $x_0$, $x_1$, $\ldots$, $x_5$, $y_0$, $y_1$, $\ldots$, $y_5$.) We define our $A_3$ and $B_3$ by $\alpha_3A_2 - \alpha_2A_3 = V$ and $\alpha_2B_3 - \alpha_3B_2 = U$. Direct computation shows that both of them are indeed elements of~$\mathcal{A}[\mathbf{x}, \mathbf{y}]$. \end{myproof}

(The question of how one might go about finding such $A$'s and $B$'s, given $U$, $V$, and $D$, will be addressed in Section \ref{o7}.)

This is the load-bearing component of our proof for Theorems \ref{o6r} and \ref{o6t}. The rest of the argument can be sketched, in very broad strokes, as follows: Each contiguous diamond minor of size $4 \times 4$ in $s \times t$ is generically nonzero. This fact and Lemma \ref{o6cm} together set up an application of Lemma \ref{cm}. We do apply it, and we arrive at Theorem \ref{o6t}. Finally, we set $s = t$, and we obtain Theorem \ref{o6r} as well.

The remainder of this section will be devoted to filling in all of the technical details which are missing from this sketch. For this purpose, we are going to employ the tools of Section \ref{algp}.

First we are going to develop a more constructive understanding of the admissible $(\mathbf{a}, \mathbf{s}, \mathbf{t})$. We temporarily set aside $s_0$ and $s_5$ as well as the indeterminates $x_0$ and $x_5$ associated with them, and we denote $\mathbf{s}^\star = (s_1, s_2, s_3, s_4)$ as well as $\mathbf{x}^\star = (x_1, x_2, x_3, x_4)$.

\begin{lemma} \label{o6g} For a generic $(\mathbf{a}, \mathbf{s}^\star, \mathbf{t}) \in \mathbb{C}^{13}$, there are exactly eight choices of $s_0$ and $s_5$ which yield an admissible $(\mathbf{a}, \mathbf{s}, \mathbf{t})$. Furthermore, exactly six out of these eight choices satisfy $s_0 \neq 0$ and $s_5 \neq 0$. \end{lemma} 

We prepare for the proof by setting up a good deal of notation. Let $R = \res_{x_5}(U, V)$. Direct computation shows that $R$ factors as $x_0 \cdot (\alpha_1x_0x_4 + \alpha_2x_1x_3 + \alpha_3x_2^2) \cdot R_\divideontimes$. We collect terms in $x_0$ to get $R = R_1x_0 + R_2x_0^2 + \cdots + R_8x_0^8$. We also collect terms in $x_5$ to get $U = U_0 + U_1x_5 + U_2x_5^2$ and $V = V_0 + V_1x_5 + V_2x_5^2$ as well as $W = V_2U - U_2V = W_0 + W_1x_5$.

Given any $P \in \mathcal{A}[\mathbf{x}, \mathbf{y}]$, we write $\widehat{P}$ for the polynomial $P(\mathbf{a}, x_0, \mathbf{s}^\star, x_5, \mathbf{t})$ formed by plugging the components of $(\mathbf{a}, \mathbf{s}^\star, \mathbf{t})$ into the corresponding indeterminates of $P$. Notice that the images under this transformation of $R$, $U_0$, $U_1$, $U_2$, $V_0$, $V_1$, $V_2$, $W_0$, $W_1$ are all in $\mathbb{C}[x_0]$.

\begin{myproof} For a generic $(\mathbf{a}, \mathbf{s}^\star, \mathbf{t}) \in \mathbb{C}^{13}$, we get that:

(i) $\widehat{R}$ is of full degree, as $R_8 \neq \zp$;

(ii) The roots of $\widehat{R}$ are all of unit multiplicity, as $\disc_{x_0}(R) \neq \zp$;

(iii) $\widehat{R}$ and $\widehat{V_2}$ do not share any common roots, as $\res_{x_0}(R, V_2) \neq \zp$;

(iv) For each root $s_0$ of $\widehat{R}$, the polynomials $\widehat{U}(s_0, x_5)$ and $\widehat{V}(s_0, x_5)$ share a common root -- as the latter is of full degree by (iii), and $R = \res_{x_5}(U, V)$;

(v) $\widehat{R}$ and $\widehat{W_1}$ do not share any common roots, as $\res_{x_0}(R, W_1) \neq \zp$;

(vi) For each root $s_0$ of $\widehat{R}$, the polynomials $\widehat{U}(s_0, x_5)$ and $\widehat{V}(s_0, x_5)$ share a unique common root given by $s_5 = -\widehat{W_0}(s_0)/\widehat{W_1}(s_0)$; indeed, the existence of a common root is guaranteed by (iv), while (v) ensures the formula's validity.

This analysis confirms the first part of Lemma \ref{o6g}. (It also tells us how, given $(\mathbf{a}, \mathbf{s}^\star, \mathbf{t})$, to find all $s_0$ and $s_5$ which lift it into an admissible $(\mathbf{a}, \mathbf{s}, \mathbf{t})$.) We continue on to the second part.

Of course, the factor $x_0$ of $R$ yields $s_0 = 0$. Similarly, the factor $\alpha_1x_0x_4 + \alpha_2x_1x_3 + \alpha_3x_2^2$ yields $s_0 = -(a_2s_1s_3 + a_3s_2^2)/a_1s_4$ and $s_5 = -\widehat{W_0}(s_0)/\widehat{W_1}(s_0) = 0$ in the generic setting of (i)--(vi). On the other hand, for a generic $(\mathbf{a}, \mathbf{s}^\star, \mathbf{t}) \in \mathbb{C}^{13}$, all roots of $\widehat{R_\divideontimes}$ yield nonzero $s_0$ and $s_5$ since $R_1 \neq \zp$ and $\res_{x_0}(R_\divideontimes, W_0) \neq \zp$. \end{myproof}

We can now spell out precisely what is meant by the term ``generically'' in the statement of Theorem \ref{o6t}. Consider any property $\mathcal{P}$ which two sequences of complex numbers might or might not possess. In light of Lemma \ref{o6g}, we define ``two order-$6$ Somos twins $s$ and $t$ with nonzero terms generically satisfy $\mathcal{P}$'' to mean ``for a generic $(\mathbf{a}, \mathbf{s}^\star, \mathbf{t}) \in \mathbb{C}^{13}$ and for every choice of $s_0 \neq 0$ and $s_5 \neq 0$ which lifts it into an admissible $(\mathbf{a}, \mathbf{s}, \mathbf{t})$, it holds that all order-$6$ Somos twins $s$ and $t$ with nonzero terms which are based on $(\mathbf{a}, \mathbf{s}, \mathbf{t})$ satisfy $\mathcal{P}$''.

From now on, we will be working with the system $\{R_\divideontimes, W\}$ in place of $\{U, V\}$. The point of this is simply to filter out the admissible $(\mathbf{a}, \mathbf{s}, \mathbf{t})$ where either $s_0 = 0$ or $s_5 = 0$.

Our next order of business will be to verify that certain expressions are generically nonzero over the solution space of $\{R_\divideontimes, W\}$. (For the same augmented notion of genericity as above.) These expressions will be listed explicitly in Lemma \ref{o6nz} below. Before we can get there, though, we must make some preliminary remarks.

One subtlety of Lemma \ref{o6nz} is that we are making a separate genericity statement for each expression. So, in particular, each expression is associated with its own implicit system of non-degeneracy conditions. We do not claim that there exists a single implicit system of non-degeneracy conditions which would resolve the matter for all of these expressions simultaneously.

The following device will be helpful: Let $P$ be any polynomial in $\mathcal{A}[\mathbf{x}, \mathbf{y}]$. We collect terms in $x_5$ to get $P = P_0 + P_1x_5 + \cdots + P_kx_5^k$, and we set $\mu(P) = W_1^k P(-W_0/W_1) = P_0W_1^k - P_1W_0W_1^{k - 1} + \cdots \pm P_kW_0^k$. The point of this definition is to eliminate $x_5$ from $P$ in a way consistent with the system $\{R_\divideontimes, W\}$.

The main construction in the proof involves finite fields. Given a prime $p$, we write $\mathbb{F}_p$ for the finite field of size $p$. In order to avoid confusion between the integers and the elements of $\mathbb{F}_p$, in this context we write $\overline{k}$ for the residue class of the integer $k$ modulo $p$, viewed as an element~of~$\mathbb{F}_p$.

\begin{lemma} \label{o6nz} Suppose that $P$ is either the numerator of some term of $S$, or the numerator of some term of $T$, or the numerator of some contiguous diamond minor of size $4 \times 4$ in $S \times T$. Then, for a generic $(\mathbf{a}, \mathbf{s}^\star, \mathbf{t}) \in \mathbb{C}^{13}$ and for every choice of $s_0 \neq 0$ and $s_5 \neq 0$ which lifts it into an admissible $(\mathbf{a}, \mathbf{s}, \mathbf{t})$, it holds that $P(\mathbf{a}, \mathbf{s}, \mathbf{t}) \neq 0$. \end{lemma} 

\begin{myproof} We are done if $\mu(P) \neq \zp$ and $\res_{x_0}(R_\divideontimes, \mu(P)) \neq \zp$. Suppose not, for the sake of contradiction.

Let $p$ be a prime and choose $(\mathbf{a}, \mathbf{s}^\star, \mathbf{t}) \in \mathbb{F}_p^{13}$ so that $\widehat{R_\divideontimes}$ is of full degree, it splits into linear factors, and it does not share any common roots with $\widehat{W_1}$. Since by assumption either $\mu(P) = \zp$ or $\res_{x_0}(R_\divideontimes, \mu(P)) = \zp$, we get that $\widehat{R_\divideontimes}$ and $\widehat{\mu(P)}$ share a common root $s_0$. We set $s_5 = -\widehat{W_0}(s_0)/\widehat{W_1}(s_0)$, and we arrive at $P(\mathbf{a}, \mathbf{s}, \mathbf{t}) = \overline{0}$.

What remains, for a contradiction, is to exhibit a prime $p$ and a suitable $(\mathbf{a}, \mathbf{s}^\star, \mathbf{t}) \in \mathbb{F}_p^{13}$ such that every root $s_0$ of $\widehat{R_\divideontimes}$ together with its corresponding $s_5$ makes all terms of $S$ and $T$ as well as all contiguous diamond minors of size $4 \times 4$ in $S \times T$ well-defined and nonzero.

We claim that $p = 19$ and $\mathbf{a} = (\overline{1}, \overline{1}, \overline{1})$, $\mathbf{s}^\star = (\overline{1}, \overline{1}, \overline{1}, \overline{1})$, $\mathbf{t} = (\overline{1}, \overline{1}, \overline{4}, \overline{4}, \overline{1}, \overline{1})$ work. In this setting, $\widehat{R_\divideontimes} = \overline{5}(x_0 - \overline{15})^6$. So we only need to consider $s_0 = s_5 = \overline{15}$. The coefficients and seeds thus obtained yield two doubly infinite order-$6$ Somos sequences $s$ and $t$ over $\mathbb{F}_{19}$, both of which are periodic with period length $612$, and in both of which all terms are nonzero.

Because of the periodicity, we only need to examine $612^2 = 374544$ contiguous diamond minors of size $4 \times 4$ in $s \times t$ so as to confirm that all contiguous diamond minors of size $4 \times 4$ in $s \times t$ are non-vanishing. This is doable, in principle, but there is also a better way. By means of Proposition \ref{dia}, we can bring the number of contiguous diamond minors to be examined down to a mere $612$. The verification does go through. \end{myproof}

The proof of Theorem \ref{o6t} does not pose much difficulty anymore.

\begin{myproof}[Proof of Theorem \ref{o6t}] Fix any diamond sub-matrix $M$ of size $5 \times 5$ in $S \times T$, and let $D$ be the numerator of the corresponding diamond minor $\det M$. Clearly, it suffices to show that $\widehat{R_\divideontimes}$ divides $\widehat{\mu(D)}$ for a generic $(\mathbf{a}, \mathbf{s}^\star, \mathbf{t}) \in \mathbb{C}^{13}$ -- provided, though, that the implicit system of non-degeneracy conditions does not depend on $M$.

Let $e'_0 < e'_1 < \cdots < e'_4$ and $e''_0 < e''_1 < \cdots < e''_4$ be the offsets of $M$. Consider also the diamond sub-matrix $M_\text{Hull}$ in $S \times T$ with offsets $e'_0$, $e'_0 + 2$, $e'_0 + 4$, $\ldots$, $e'_4$ and $e''_0$, $e''_0 + 2$, $e''_0 + 4$, $\ldots$, $e''_4$. For a generic $(\mathbf{a}, \mathbf{s}^\star, \mathbf{t}) \in \mathbb{C}^{13}$ and for every choice of $s_0 \neq 0$ and $s_5 \neq 0$ which lifts it into an admissible $(\mathbf{a}, \mathbf{s}, \mathbf{t})$, we get that:

(i) By Lemma \ref{o6nz}, each entry of $M_\text{Hull}(\mathbf{a}, \mathbf{s}, \mathbf{t})$ is well-defined and each contiguous ordinary minor of size $4 \times 4$ in $M_\text{Hull}(\mathbf{a}, \mathbf{s}, \mathbf{t})$ is nonzero;

(ii) By Lemma \ref{o6cm}, each contiguous ordinary minor of size $5 \times 5$ in $M_\text{Hull}(\mathbf{a}, \mathbf{s}, \mathbf{t})$ vanishes;

(iii) By (i) and (ii) together with Lemma \ref{cm}, in fact all ordinary minors of size $5 \times 5$ in $M_\text{Hull}(\mathbf{a}, \mathbf{s}, \mathbf{t})$ must vanish.

It is a direct corollary of (iii) that $\widehat{R_\divideontimes}$ divides $\widehat{\mu(D)}$ for a generic $(\mathbf{a}, \mathbf{s}^\star, \mathbf{t}) \in \mathbb{C}^{13}$. However, the implicit non-degeneracy system in this corollary depends on $M$. We rectify matters as follows: By Lemmas \ref{gen} and \ref{div}, in fact $\widehat{R_\divideontimes}$ divides $\widehat{\mu(D)}$ for all $(\mathbf{a}, \mathbf{s}^\star, \mathbf{t}) \in \mathbb{C}^{13}$ such that $\widehat{R_\divideontimes} \neq \zp$. This does give us an implicit system of non-degeneracy conditions fully independent of $M$. \end{myproof}

By the proofs of Lemma \ref{o6g} and Theorem \ref{o6t}, we obtain that one implicit system of non-degeneracy conditions for Theorem \ref{o6t} is given by $\Gamma = \{R_1, R_8, \res_{x_0}(R_\divideontimes, W_0), \res_{x_0}(R_\divideontimes, W_1)\}$.

We proceed, finally, to derive Theorem \ref{o6r} from Theorem \ref{o6t}.

\begin{myproof}[Proof of Theorem \ref{o6r}] Fix any diamond minor $\Delta$ of size $5 \times 5$ in $S \times S$. We wish to show that $\Delta = \zp$. This will follow if we can demonstrate that $\Delta(\mathbf{a}, \mathbf{s}) = 0$ for a generic $(\mathbf{a}, \mathbf{s}) \in \mathbb{C}^9$.

Consider the diamond minor $\Delta_\text{Twin}$ in $S \times T$ with the same offsets as $\Delta$. Then $\Delta(\mathbf{a}, \mathbf{s}) = \Delta_\text{Twin}(\mathbf{a}, \mathbf{s}, \mathbf{s})$ so long as both sides are well-defined -- which they are, generically. On the other hand, $\Delta_\text{Twin}(\mathbf{a}, \mathbf{s}, \mathbf{s}) = 0$ whenever the left-hand side is well-defined and $(\mathbf{a}, \mathbf{s}^\star, \mathbf{s})$ is non-degenerate in the sense of the implicit system of non-degeneracy conditions $\Gamma$. What remains is to verify that, for a generic $(\mathbf{a}, \mathbf{s}) \in \mathbb{C}^9$, the point $(\mathbf{a}, \mathbf{s}^\star, \mathbf{s}) \in \mathbb{C}^{13}$ is indeed non-degenerate in this manner. Or, equivalently, that $P(\boldsymbol{\alpha}, \mathbf{x}^\star, \mathbf{x}) \neq \zp$ for all $P \in \Gamma$. This is straightforward. \end{myproof}

\section{Order 7} \label{o7}

Here, we prove Theorems \ref{o7r} and \ref{o7t}. The argument will follow very closely in the footsteps of the order-$6$ one we presented in Section \ref{o6}. For this reason, we are going to focus most of all on the small number of differences.

Let $s$ and $t$ be two order-$7$ Somos sequences, both with coefficients $\mathbf{a}$, and with seeds $\mathbf{s}$ and $\mathbf{t}$, respectively. We clone the order-$7$ master Somos sequence $S$ and the order-$7$ Somos invariants obtained in Section \ref{inv} again, using the new formal indeterminates $y_0$, $y_1$, $\ldots$, $y_6$. Let $U = \Pi_Y\Phi_X - \Pi_X\Phi_Y$ be the numerator of $F_X - F_Y$, and similarly let $V = \Pi_Y\Psi_X - \Pi_X\Psi_Y$ be the numerator of $G_X - G_Y$. We define an admissible $(\mathbf{a}, \mathbf{s}, \mathbf{t}) \in \mathbb{C}^{17}$ as before.

\begin{lemma} \label{o7cm} Suppose that $s$ and $t$ are twins with nonzero terms for which $a_1$ and $a_3$ are not both zero. Then all contiguous half-diamond minors of size $5 \times 5$ in the matrix $s \times t$ vanish. \end{lemma} 

The condition that $(a_1, a_3) \neq (0, 0)$ cannot be omitted; indeed, setting $(\alpha_1, \alpha_3) = (0, 0)$ in the argument below causes $U$ and $V$ to vanish but not $D$.

\begin{myproof} Let $\Delta$ be the contiguous half-diamond minor of size $5 \times 5$ in $S \times T$ with offsets $(-4, -2, 0, 2, 4)$ and $(-2, 2, 6, 10, 14)$. Let also $\Delta = D/x_0x_6y_0y_6\Pi_X\Pi_Y$. (The denominator of $\Delta$ in lowest terms is in fact $x_0^2x_1x_2x_4x_5x_6^2y_0^2y_1y_2y_4y_5y_6^2$. The numerator $D$ is homogeneous of degree $28$, with $3989$ summands.)

We claim that there exist $A_1$ and $B_1$ in $\mathcal{A}[\mathbf{x}, \mathbf{y}]$ with $A_1U + B_1V = \alpha_1D$; as well as $A_3$ and $B_3$ in $\mathcal{A}[\mathbf{x}, \mathbf{y}]$ with $A_3U + B_3V = \alpha_3D$. Clearly, that would suffice.

\begin{table} \footnotesize \centering \begin{tabular}{cccccc} $00232002003002$ & $\alpha_1^2\alpha_3^2$ & $11022010032110$ & $\alpha_1\alpha_2\alpha_3^2$ & $11120020211201$ & $\alpha_1^3\alpha_3$\\ $00321102003002$ & $\alpha_1^2\alpha_2\alpha_3$ & $11022010041020$ & $\alpha_1\alpha_2^2\alpha_3$ & $11120020300301$ & $\alpha_1^3\alpha_2$\\ $01123002003002$ & $\alpha_1^2\alpha_2\alpha_3$ & $11022010121210$ & $\alpha_1\alpha_2^2\alpha_3$ & $11120021013101$ & $\alpha_1^3\alpha_3$\\ $01212102003002$ & $\alpha_1^2\alpha_2^2$ & $11022010122101$ & $\alpha_1^2\alpha_3^2$ & $11120021022011$ & $\alpha_1^3\alpha_2$\\ $01221012003002$ & $\alpha_1^3\alpha_3$ & $11022010130120$ & $\alpha_1\alpha_2^3$ & $11120021101310$ & $\alpha_1^3\alpha_3$\\ $02112011111111$ & $\alpha_2^2\alpha_3$ & $11022010131011$ & $2\alpha_1^2\alpha_2\alpha_3$ & $11120021110220$ & $\alpha_1^3\alpha_2$\\ $02112012003002$ & $\alpha_1^3\alpha_2$ & $11022010220111$ & $\alpha_1^2\alpha_2^2$ & $11120021111111$ & $\alpha_1\alpha_2^2$\\ $02112012011021$ & $\alpha_1\alpha_2\alpha_3$ & $11022010221002$ & $\alpha_1^3\alpha_3$ & $11120021201102$ & $\alpha_1^2\alpha_2$\\ $02201111111111$ & $\alpha_2^3$ & $11022011013101$ & $\alpha_1^2\alpha_2\alpha_3$ & $11120022003002$ & $\alpha_1^2\alpha_2$\\ $02201112011021$ & $\alpha_1\alpha_2^2$ & $11022011021120$ & $\alpha_1^2\alpha_2\alpha_3$ & $11200211111111$ & $\alpha_1\alpha_2^2$\\ $03101021111111$ & $\alpha_1\alpha_2^2$ & $11022011030030$ & $\alpha_1^2\alpha_2^2$ & $11200212011021$ & $\alpha_1^2\alpha_2$\\ $03101022011021$ & $\alpha_1^2\alpha_2$ & $11022011111111$ & $\alpha_2^3$ & $12002111111111$ & $\alpha_1\alpha_2^2$\\ $10122102003002$ & $\alpha_1^3\alpha_3$ & $11022012003002$ & $\alpha_1\alpha_2^2$ & $12002111201102$ & $\alpha_1^2\alpha_2$\\ $10131012003002$ & $\alpha_1\alpha_2\alpha_3$ & $11022012011021$ & $\alpha_1\alpha_2^2$ & $20021110032110$ & $\alpha_1^2\alpha_3^2$\\ $10211201111111$ & $\alpha_2^2\alpha_3$ & $11102201111111$ & $\alpha_2^3$ & $20021110041020$ & $\alpha_1^2\alpha_2\alpha_3$\\ $10211201201102$ & $\alpha_1\alpha_2\alpha_3$ & $11102201201102$ & $\alpha_1\alpha_2^2$ & $20021110121210$ & $\alpha_1^2\alpha_2\alpha_3$\\ $10211202003002$ & $\alpha_1^3\alpha_2$ & $11111110122101$ & $\alpha_1^2\alpha_2\alpha_3$ & $20021110130120$ & $\alpha_1^2\alpha_2^2$\\ $10220110112300$ & $\alpha_1\alpha_2\alpha_3^2$ & $11111110131011$ & $\alpha_1^2\alpha_2^2$ & $20021110131011$ & $\alpha_1^3\alpha_3$\\ $10220110121210$ & $\alpha_1\alpha_2^2\alpha_3$ & $11111110221002$ & $\alpha_1^3\alpha_2$ & $20021110220111$ & $\alpha_1^3\alpha_2$\\ $10220110201400$ & $\alpha_1\alpha_2^2\alpha_3$ & $11111111012210$ & $\alpha_1^2\alpha_2\alpha_3$ & $20021111013101$ & $\alpha_1^3\alpha_3$\\ $10220110210310$ & $\alpha_1\alpha_2^3$ & $11111111013101$ & $\alpha_1^2\alpha_2^2$ & $20021111021120$ & $\alpha_1^3\alpha_3$\\ $10220110211201$ & $\alpha_1^2\alpha_2\alpha_3$ & $11111111022011$ & $\alpha_1^3\alpha_3$ & $20021111030030$ & $\alpha_1^3\alpha_2$\\ $10220110300301$ & $\alpha_1^2\alpha_2^2$ & $11111111101310$ & $\alpha_1^2\alpha_2^2$ & $20021111102201$ & $\alpha_1^3\alpha_2$\\ $10220111012210$ & $\alpha_1^2\alpha_3^2$ & $11111111102201$ & $\alpha_1^3\alpha_3$ & $20021111111111$ & $\alpha_1\alpha_2^2$\\ $10220111013101$ & $\alpha_1^2\alpha_2\alpha_3$ & $11111111201102$ & $2\alpha_1^2\alpha_3$ & $20021112003002$ & $\alpha_1^2\alpha_2$\\ $10220111101310$ & $2\alpha_1^2\alpha_2\alpha_3$ & $11111112001220$ & $\alpha_1^3\alpha_2$ & $20021112011021$ & $\alpha_1^2\alpha_2$\\ $10220111110220$ & $\alpha_1^2\alpha_2^2$ & $11111112011021$ & $2\alpha_1^2\alpha_3$ & $20101301111111$ & $\alpha_1\alpha_2^2$\\ $10220111111111$ & $\alpha_2^3$ & $11120020112300$ & $\alpha_1^2\alpha_3^2$ & $20101301201102$ & $\alpha_1^2\alpha_2$\\ $10220111201102$ & $\alpha_1\alpha_2^2$ & $11120020121210$ & $\alpha_1^2\alpha_2\alpha_3$ & $21010121111111$ & $\alpha_2^2$\\ $10220112001220$ & $\alpha_1^3\alpha_3$ & $11120020201400$ & $\alpha_1^2\alpha_2\alpha_3$ & $21010121201102$ & $\alpha_1\alpha_2$\\ $10220112003002$ & $\alpha_1\alpha_2^2$ & $11120020210310$ & $\alpha_1^2\alpha_2^2$ & $21010122011021$ & $\alpha_1\alpha_2$ \end{tabular} \caption{} \label{o7a1} \end{table}

\begin{table} \footnotesize \centering \begin{tabular}{cccccc} $01221012002111$ & $\alpha_1\alpha_3^2$ & $10131011111111$ & $\alpha_1^3\alpha_3$ & $11111112100121$ & $\alpha_2^2$\\ $01310112002111$ & $\alpha_1\alpha_2\alpha_3$ & $11013101112002$ & $\alpha_1\alpha_2\alpha_3$ & $11120021022011$ & $\alpha_1\alpha_2^2$\\ $02210022002111$ & $\alpha_1^2\alpha_3$ & $11111110122101$ & $\alpha_2^2\alpha_3$ & $11200211102201$ & $\alpha_1\alpha_2^2$\\ $10122101112002$ & $\alpha_1\alpha_3^2$ & $11111110131011$ & $\alpha_2^3$ & $11200212002111$ & $\alpha_1^2\alpha_2$\\ $10131010023200$ & $\alpha_1\alpha_3^3$ & $11111110221002$ & $\alpha_1\alpha_2^2$ & $12002111022011$ & $\alpha_1\alpha_2^2$\\ $10131010032110$ & $\alpha_1\alpha_2\alpha_3^2$ & $11111111012210$ & $\alpha_2^2\alpha_3$ & $12002111112002$ & $\alpha_1^2\alpha_2$\\ $10131010112300$ & $\alpha_1\alpha_2\alpha_3^2$ & $11111111013101$ & $\alpha_2^3$ & $12011022100121$ & $\alpha_1\alpha_2$\\ $10131010121210$ & $\alpha_1\alpha_2^2\alpha_3$ & $11111111022011$ & $\alpha_1\alpha_2\alpha_3$ & $20012201112002$ & $\alpha_1^2\alpha_3$\\ $10131010122101$ & $\alpha_1^2\alpha_3^2$ & $11111111101310$ & $\alpha_2^3$ & $20021111102201$ & $\alpha_1\alpha_2^2$\\ $10131010211201$ & $\alpha_1^2\alpha_2\alpha_3$ & $11111111102201$ & $\alpha_1\alpha_2\alpha_3$ & $20030021111111$ & $\alpha_1^2\alpha_2$\\ $10131011012210$ & $\alpha_1^2\alpha_3^2$ & $11111111210012$ & $\alpha_2^2$ & $20110211210012$ & $\alpha_1\alpha_2$\\ $10131011021120$ & $\alpha_1^2\alpha_2\alpha_3$ & $11111112001220$ & $\alpha_1\alpha_2^2$ & $21010121111111$ & $\alpha_1\alpha_3$ \end{tabular} \caption{} \label{o7b1} \end{table}

We compress these $A$'s and $B$'s as in the proof of Lemma \ref{o6cm}. Our $A_1$ and $B_1$ are the skew-symmetrisations of the polynomials shown in Tables \ref{o7a1} and \ref{o7b1}, respectively. We define our $A_3$ and $B_3$ by $\alpha_3A_1 - \alpha_1A_3 = \alpha_2^2V$ and $\alpha_1B_3 - \alpha_3B_1 = \alpha_2^2U$. Direct computation shows that both of them are indeed elements of $\mathcal{A}[\mathbf{x}, \mathbf{y}]$. \end{myproof}

How does one find such $A$'s and $B$'s, given $U$, $V$, and $D$? For polynomials of the size we are dealing with, this is not a trivial problem, and we go on a brief digression to address it.

Let $\mathbb{P}$ be a polynomial ring, and fix $P_1$, $P_2$, $\ldots$, $P_k \in \mathbb{P}$. The set $\mathbb{I}$ of all $Q = C_1P_1 + C_2P_2 + \cdots + C_kP_k$, with $C_1$, $C_2$, $\ldots$, $C_k \in \mathbb{P}$, is known as the \emph{ideal} generated by $P_1$, $P_2$, $\ldots$, $P_k$ in $\mathbb{P}$. We say that $C_1$, $C_2$, $\ldots$, $C_k$ \emph{certify} the membership of $Q$ in $\mathbb{I}$.

The task of deciding and certifying ideal membership comes up often enough in applications that many computer algebra systems include tools for its automation. For example, the certificates in the proof of Lemma \ref{o7cm} were obtained with the help of the \texttt{lift} command of the computer algebra system \textsc{Singular} \cite{DGPS}.

Our polynomials are so big that the computation is not actually feasible directly. Still, we can make do by the following trick: Plug $x_2 = x_3 = x_4 = 1$ and $y_2 = y_3 = y_4 = 1$ into $U$, $V$, and $D$. For the new polynomials thus obtained, the computation goes through in a matter of seconds on modern hardware. We can now recover the missing exponents of $x_2$, $x_3$, $x_4$, $y_2$, $y_3$, $y_4$ in the certificates by exploiting certain linear constraints on the exponent tuples of $U$, $V$, and $D$ associated with $\mathcal{E}$, similarly to how we employed $\mathcal{E}$ in our definition of $\Upsilon_\boxtimes$ and $\Omega_\boxtimes$ in Section~\ref{inv}.

This concludes our digression on the topic of computation, and we return to our main narrative. We temporarily set aside $s_0$ and $s_6$, and we denote $\mathbf{s}^\star = (s_1, s_2, \ldots, s_5)$.

\begin{lemma} \label{o7g} For a generic $(\mathbf{a}, \mathbf{s}^\star, \mathbf{t}) \in \mathbb{C}^{15}$, there are exactly eight choices of $s_0$ and $s_6$ which yield an admissible $(\mathbf{a}, \mathbf{s}, \mathbf{t})$. Furthermore, exactly six out of these eight choices satisfy $s_0 \neq 0$ and $s_6 \neq 0$. \end{lemma} 

The proof is analogous to that of Lemma \ref{o6g}. This time around, $R = \res_{x_6}(U, V)$ factors as $x_0 \cdot (\alpha_1x_0x_5 + \alpha_2x_1x_4 + \alpha_3x_2x_3) \cdot R_\divideontimes$.

Just as in Section \ref{o6}, it is on the basis of Lemma \ref{o7g} that we augment out notion of genericity in the context of Theorem \ref{o7t}. We define ``two order-$7$ Somos twins $s$ and $t$ with nonzero terms generically satisfy $\mathcal{P}$'' to mean ``for a generic $(\mathbf{a}, \mathbf{s}^\star, \mathbf{t}) \in \mathbb{C}^{15}$ and for every choice of $s_0 \neq 0$ and $s_6 \neq 0$ which lifts it into an admissible $(\mathbf{a}, \mathbf{s}, \mathbf{t})$, it holds that all order-$7$ Somos twins $s$ and $t$ with nonzero terms which are based on $(\mathbf{a}, \mathbf{s}, \mathbf{t})$ satisfy $\mathcal{P}$''.

\begin{lemma} \label{o7nz} Suppose that $P$ is either the numerator of some term of $S$, or the numerator of some term of $T$, or the numerator of some contiguous half-diamond minor of size $4 \times 4$ in $S \times T$. Then, for a generic $(\mathbf{a}, \mathbf{s}^\star, \mathbf{t}) \in \mathbb{C}^{15}$ and for every choice of $s_0 \neq 0$ and $s_6 \neq 0$ which lifts it into an admissible $(\mathbf{a}, \mathbf{s}, \mathbf{t})$, it holds that $P(\mathbf{a}, \mathbf{s}, \mathbf{t}) \neq 0$. \end{lemma} 

\begin{myproof} Just as with Lemma \ref{o6nz}, it suffices to find a suitable finite-field construction.

We claim that $p = 29$ and $\mathbf{a} = (\overline{1}, \overline{1}, \overline{1})$, $\mathbf{s}^\star = (\overline{1}, \overline{1}, \overline{6}, \overline{1}, \overline{1})$, $\mathbf{t} = (\overline{1}, \overline{1}, \overline{2}, \overline{1}, \overline{9}, \overline{1}, \overline{1})$ work. In this setting, $\widehat{R_\divideontimes} = (x_0 - \overline{3})^2(x_0 - \overline{4})^2(x_0 - \overline{5})^2$. Our options for $(s_0, s_6)$ are, accordingly, $(\overline{3}, \overline{4})$, $(\overline{4}, \overline{3})$, $(\overline{5}, \overline{5})$. Each one of the three seeds $\mathbf{s}$ thus obtained yields a doubly infinite order-$7$ Somos sequences over $\mathbb{F}_{29}$ with nonzero terms, and so does the seed $\mathbf{t}$ as well. The period lengths of these sequences are $112$, $112$, $16$, $16$, respectively.

There are three options for the matrix $s \times t$ to consider. Because of the periodicity, over all three options we must verify the non-vanishing of only $2 \cdot 16 \cdot 112 + 2 \cdot16 \cdot 112 + 2 \cdot 16 \cdot 16 = 7680$ contiguous half-diamond minors of size $4 \times 4$ altogether. This is not a particularly difficult computation even without the optimisations enabled by Proposition \ref{dia}, and the verification does go through. \end{myproof}

The proof of Theorem \ref{o7t} is analogous to that of Theorem \ref{o6t}. The only difference worth remarking upon is that we must be a little more careful as we form $M_\text{Hull}$. (Fix any half-diamond sub-matrix $M$ of size $5 \times 5$ in $S \times T$, with offsets $e'$ and $e''$. Suppose, for concreteness, that the terms of $e'$ are pairwise congruent modulo $4$. Then we let $M_\text{Hull}$ be the half-diamond sub-matrix of $S \times T$ with offsets $e'_0$, $e'_0 + 4$, $e'_0 + 8$, $\ldots$, $e'_4$ and $e''_0$, $e''_0 + 2$, $e''_0 + 4$, $\ldots$, $e''_4$.) On the other hand, the proof of Theorem \ref{o7r} is fully analogous to the proof of Theorem \ref{o6r}.

\section{Integrality} \label{int}

The early popularity of Somos sequences owes a lot to their integrality properties. Even though the Somos recurrence involves division, the unit Somos sequence of order $n$ turns out to consist entirely of integers for all $2 \le n \le 7$. This property is preserved when we retain the unit seed but allow the coefficients to be arbitrary integers. A stronger generalisation \cite{G, FZ} dispenses with concrete numbers altogether and is phrased instead in terms of the corresponding master Somos sequences. Theorems \ref{o6i} and \ref{o7i} below state this generalisation in the setting of orders $6$~and~$7$, respectively.

A rational function whose denominator is the product of several indeterminates is known as a \emph{Laurent polynomial}. These objects may also be viewed as ``generalised polynomials'' where negative exponents are allowed. For example, $(x^2 + y^2)/xy = xy^{-1} + x^{-1}y$.

\begin{theorem} \label{o6i} Every term of the master Somos sequence of order $6$ is a Laurent polynomial. \end{theorem} 

\begin{theorem} \label{o7i} Every term of the master Somos sequence of order $7$ is a Laurent polynomial. \end{theorem} 

Our purpose here will be to give a proof of Theorems \ref{o6i} and \ref{o7i} based exclusively on the finite-rank properties of Somos sequences. The utility of these properties in the study of arithmetical questions has already been demonstrated in \cite{U24}.

Let $\mathbb{U}$ be a unique factorisation domain and let $\mathbb{U}_\text{Frac}$ be the fraction field of $\mathbb{U}$. For example, if $\mathbb{U} = \mathbb{Z}$ then $\mathbb{U}_\text{Frac} = \mathbb{Q}$. The special case which matters for Theorems \ref{o6i} and \ref{o7i} is when $\mathbb{U} = \mathcal{A}[\mathbf{x}]$ and $\mathbb{U}_\text{Frac} = \mathcal{A}(\mathbf{x})$.

Let $u$ be a sequence in $\mathbb{U}_\text{Frac}$ indexed by $\mathbb{Z}$. (We assume that $u$ is doubly infinite purely for convenience. It is straightforward to adapt the argument to finite and singly-infinite sequences.) Let also $\Theta$ be the set of all irreducibles in $\mathbb{U}$ which divide the denominator of some term of $u$. The gist of what follows is that the matrix $u \times u$ being of finite diamond rank tells us a lot about~$\Theta$.

Let $\Theta(i)$ be the set of all irreducibles in $\mathbb{U}$ which divide the denominator of some term of the subsequence of $u$ indexed by $[-i; i]$. So, in particular, $\Theta(0) \subseteq \Theta(1) \subseteq \Theta(2) \subseteq \cdots$ and $\Theta = \bigcup_i \Theta(i)$. Let also $\Lambda$ be the set of all irreducibles in $\mathbb{U}$ which divide either the numerator of $u_0$, or both of the numerators of $u_{-1}$ and $u_1$.

Fix a positive integer $r$. We write $\Xi(i)$ for the set of all irreducibles in $\mathbb{U}$ which divide the numerator of every diamond minor of size $r \times r$ in $u \times u$ all of whose offsets are in the interval $[-i; i]$. So, in particular, $\Xi(0) \supseteq \Xi(1) \supseteq \Xi(2) \supseteq \cdots$; the nesting is in the opposite direction relative to $\Theta(i)$ because the definition of $\Theta(i)$ uses an existential quantifier, while in the definition of $\Xi(i)$ the quantifier is universal instead.

\begin{lemma} \label{df} Suppose that the terms of $u$ are nonzero and the matrix $u \times u$ is of diamond rank~$r$. Then $\Theta$ is finite. Furthermore, for each positive integer $k$, it holds that $\Theta \subseteq \Lambda \cup \Xi(k) \cup \Theta(k + 1)$. \end{lemma} 

\begin{myproof} The first part is a corollary of the second one, and so we focus on the latter.

Fix a positive integer $k$. We will show, by induction on $i$, that $\Theta(i) \subseteq \Lambda \cup \Xi(k) \cup \Theta(k + 1)$ for all $i$. This is clear when $0 \le i \le k + 1$. Suppose now that $i \ge k + 2$, and let $p$ be any irreducible in $\mathbb{U}$ which divides the denominator of $u_i$ (the case of $u_{-i}$ is analogous) but does not belong to $\Theta(i - 1)$ or $\Lambda$. We wish to verify that $p$ belongs to $\Xi(k)$.

Consider any diamond sub-matrix $M$ of size $r \times r$ in $u \times u$ all of whose offsets are in the interval $[-k; k]$. Consider also any position $(i, j)$ in $u \times u$, with $j \in \{-1, 0, 1\}$, such that its offsets are of the same parity as the offsets of $M$. We adjoin the diagonal and the anti-diagonal through $(i, j)$ to the diagonals and the anti-diagonals of $M$, and we obtain a diamond sub-matrix $M_\text{Ext}$ of size $(r + 1) \times (r + 1)$ in $u \times u$.

Since $u \times u$ is of diamond rank $r$, we get that $\det M_\text{Ext}$ vanishes. Expanding, we get that $u_iu_j \cdot \det M$ is a signed sum of $(r + 1)! - r!$ products all of whose multiplicands are among $u_{1 - i}$, $u_{2 - i}$, $\ldots$, $u_{i - 2}$, $u_{i - 1}$. Since $p \not \in \Theta(i - 1)$, it follows that $p$ divides the numerator of $u_j \cdot \det M$ for all $j \in \{-1, 0, 1\}$ of the right parity. Regardless of which parity it is, by virtue of $p \not \in \Lambda$ we obtain that $p$ must divide the numerator of $\det M$. Since the same reasoning applies to all diamond sub-matrices $M$ of size $r \times r$ in $u \times u$ all of whose offsets are in the interval $[-k; k]$, we conclude that $p$ does indeed belong to $\Xi(k)$, as desired. This completes the induction step. \end{myproof}

Theorem \ref{o6i} now boils down to a moderate amount of computation:

\begin{myproof}[Proof of Theorem \ref{o6i}] We apply Lemma \ref{df} with $\mathbb{U} = \mathcal{A}[\mathbf{x}]$, $\mathbb{U}_\text{Frac} = \mathcal{A}(\mathbf{x})$, $r = 4$, $S$ being the order-$6$ master Somos sequence, and $u_i = S_{i + 2}$. Let $k = 7$. Direct computation shows that $\Lambda = \{x_2\}$ and $\Theta(8) = \{x_0, x_1, \ldots, x_5\}$. Let $e = (1, 3, 5, 7)$, $e' = (-3, -1, 1, 3)$, $e'' = (-1, 1, 3, 5)$. Consider the diamond minor $\Delta'$ of $u \times u$ with offsets $e$ and $e'$, as well as the diamond minor $\Delta''$ of $u \times u$ with offsets $e$ and $e''$. The numerator of each one of $\Delta'$ and $\Delta''$ is a homogeneous polynomial of degree $17$ with $197$ summands. Direct computation shows that these two numerators are relatively prime, and so $\Xi(7) = \varnothing$. Thus $\Theta = \Theta(8)$. \end{myproof}

We proceed next to develop a half-diamond analogue of Lemma \ref{df}. Let $\Lambda_\text{Half}$ be the set of all irreducibles in $\mathbb{U}$ which divide either the numerator of $u_{-1}$, or that of $u_0$, or that of $u_1$, or both of the numerators of $u_{-2}$ and $u_2$. Let also $\Xi_\text{Half}(i)$ be the set of all irreducibles in $\mathbb{U}$ which divide the numerator of every half-diamond minor of size $r \times r$ in $u \times u$ all of whose offsets are in the interval~$[-i; i]$.

\begin{lemma} \label{hdf} Suppose that the terms of $u$ are nonzero and the matrix $u \times u$ is of half-diamond rank $r$. Then $\Theta$ is finite. Furthermore, for each positive integer $k$, it holds that $\Theta \subseteq \Lambda_\textnormal{Half} \cup \Xi_\textnormal{Half}(k) \cup \Theta(k + 2)$. \end{lemma} 

\begin{myproof} We argue as with Lemma \ref{df}. We only highlight the differences which have to do with the minors being half-diamond instead of diamond ones. Suppose, for concreteness, that the diagonal offsets of $M$ are pairwise congruent modulo $4$. Consider any position $(i, j)$ in $u \times u$, with $j \in \{-2, -1, 0, 1, 2\}$, for which the offset of the diagonal through it is congruent modulo $4$ to the offsets of the diagonals of $M$. By adjoining the diagonal and the anti-diagonal through $(i, j)$ to the diagonals and the anti-diagonals of $M$, we obtain a half-diamond sub-matrix $M_\text{Ext}$ of size $(r + 1) \times (r + 1)$ in $u \times u$. The rest of the argument goes as before. \end{myproof}

Once again, Lemma \ref{hdf} boils Theorem \ref{o7i} down to a moderate amount of computation:

\begin{myproof}[Proof of Theorem \ref{o7i}] We apply Lemma \ref{hdf} with $\mathbb{U} = \mathcal{A}[\mathbf{x}]$, $\mathbb{U}_\text{Frac} = \mathcal{A}(\mathbf{x})$, $r = 4$, $S$ being the order-$7$ master Somos sequence, and $u_i = S_{i + 3}$. Let $k = 8$. Direct computation shows that $\Lambda_\text{Half} = \{x_2, x_3, x_4\}$ and $\Theta(10) = \{x_0, x_1, \ldots, x_6\}$. Let $e = (0, 2, 4, 6)$, $e' = (-8, -4, 0, 4)$, $e'' = (-4, 0, 4, 8)$. Consider the half-diamond minor $\Delta'$ of $u \times u$ with offsets $e$ and $e'$, as well as the half-diamond minor $\Delta''$ of $u \times u$ with offsets $e$ and $e''$. The numerator of each one of $\Delta'$ and $\Delta''$ is a homogeneous polynomial of degree $17$ with $191$ summands. Direct computation shows that these two numerators are relatively prime, and so $\Xi_\text{Half}(8) = \varnothing$. Thus $\Theta = \Theta(10)$. \end{myproof}

\section{Lower Orders} \label{low}

The lower orders $2$, $3$, $4$, $5$ exhibit decimation properties similar to the ones of Propositions \ref{o6e}~and~\ref{o7e} as well as finite-rank properties similar to the ones of Theorems \ref{o6r}--\ref{o7t}. We focus on the latter as it is not too difficult to derive the former from them.

We begin with orders $2$ and $3$. The following analogues are clear from the definition of a Somos sequence:

\begin{theorem} Let $s$ and $t$ be two Somos sequences of order $2$, with the same coefficient and with nonzero terms. Then the matrix $s \times t$ is of unit diamond rank. \end{theorem}

\begin{theorem} Let $s$ and $t$ be two Somos sequences of order $3$, with the same coefficient and with nonzero terms. Then the matrix $s \times t$ is of unit half-diamond rank. \end{theorem}

We go on to orders $4$ and $5$. For a more in-depth discussion of them, we refer readers to \cite{H05}, \cite{H07}, \cite{U19}, \cite{TT} as well as their bibliographies. Here, we limit ourselves to the analogues of Theorems \ref{o6r}--\ref{o7t}:

\begin{theorem} \label{o4r} Let $s$ be a Somos sequence of order $4$ with nonzero terms. Then the matrix $s \times s$ is of diamond rank at most $2$. \end{theorem}

\begin{theorem} \label{o4t} Let $s$ and $t$ be two twinned Somos sequences of order $4$ with nonzero terms. Then, generically, the matrix $s \times t$ is of diamond rank at most $2$. \end{theorem}

\begin{theorem} \label{o5r} Let $s$ be a Somos sequence of order $5$ with nonzero terms. Then the matrix $s \times s$ is of half-diamond rank at most $2$. \end{theorem}

\begin{theorem} \label{o5t} Let $s$ and $t$ be two twinned Somos sequences of order $5$ with nonzero terms. Then, generically, the matrix $s \times t$ is of half-diamond rank at most $2$. \end{theorem}

The proofs go as in Section \ref{o6}, except that both the logical structure and the necessary computations are much simpler. This is due chiefly to the fact that now there is only one invariant to deal with. We proceed to review the key changes relative to orders $6$ and $7$. We will consider orders $4$ and $5$ in parallel, as they behave very similarly.

For order $4$, we let $\Phi_X = \Phi_4$ and $\Pi_X = \Pi_4$. For order $5$, we let $\Phi_X = \Phi_5$ and $\Pi_X = \Pi_5$. Either way, we let $\Phi_Y$ and $\Pi_Y$ be the ``clones'' of $\Phi_X$ and $\Pi_X$, and we set $U = \Pi_Y\Phi_X - \Pi_X\Phi_Y$.

For the analogue of Lemma \ref{o6cm}, with order $4$ we let $D$ be the numerator of the diamond minor of size $3 \times 3$ in $S \times T$ with offsets $(-2, 0, 2)$ and $(2, 4, 6)$. For order $5$, instead $D$ becomes the half-diamond minor of size $3 \times 3$ in $S \times T$ with offsets $(-2, 0, 2)$ and $(0, 4, 8)$.

Previously, with orders $6$ and $7$, we had two invariants governing the twinning relation, and so we had to produce ideal membership certificates for $D$ with respect to the ideal generated by $U$ and $V$. Here, though, there is only one invariant, and accordingly all we need to do is verify that $U$ divides $D$. The computations can easily be carried out by hand. For order $4$, we get that $D = -U$. For order $5$, we get that $D = \alpha_2U$.

For the analogue of Lemma \ref{o6g}, we set aside $s_0$ to define $\mathbf{s}^\star = (s_1, s_2, s_3)$ with order $4$ and $\mathbf{s}^\star = (s_1, s_2, s_3, s_4)$ with order $5$. Either way, we also let $\widehat{U} = U(\mathbf{a}, x_0, \mathbf{s}^\star, \mathbf{t}) \in \mathbb{C}[x_0]$. It is straightforward to see that, for a generic $(\mathbf{a}, \mathbf{s}^\star, \mathbf{t})$, the polynomial $\widehat{U}$ has two distinct roots both of which are nonzero.

On this basis, in both Theorems \ref{o4t} and \ref{o5t}, we understand ``two Somos twins $s$ and $t$ with nonzero terms generically satisfy $\mathcal{P}$'' to mean ``for a generic $(\mathbf{a}, \mathbf{s}^\star, \mathbf{t})$ and for every choice of $s_0$ which lifts it into an admissible $(\mathbf{a}, \mathbf{s}, \mathbf{t})$, it holds that all Somos twins $s$ and $t$ with nonzero terms which are based on $(\mathbf{a}, \mathbf{s}, \mathbf{t})$ satisfy $\mathcal{P}$''.

For the analogue of Lemma \ref{o6nz}, we must exhibit finite-field constructions where all contiguous diamond or half-diamond minors of size $2 \times 2$ are well-defined and non-vanishing. For order $4$, we let $p = 11$ as well as $\mathbf{a} = (\overline{1}, \overline{1})$, $\mathbf{s}^\star = (\overline{1}, \overline{9}, \overline{1})$, $\mathbf{t} = (\overline{1}, \overline{2}, \overline{2}, \overline{1})$. Then $\widehat{U} = \overline{4}(x_0 - \overline{4})^2$. Both the seed $\mathbf{s}$ given by the double root of $\widehat{U}$ and the seed $\mathbf{t}$ yield doubly infinite order-$4$ Somos sequences over $\mathbb{F}_{11}$ with period length $5$. For order $5$, we let $p = 11$ as well as $\mathbf{a} = (\overline{1}, \overline{1})$, $\mathbf{s}^\star = (\overline{1}, \overline{1}, \overline{2}, \overline{1})$, $\mathbf{t} = (\overline{1}, \overline{1}, \overline{5}, \overline{1}, \overline{1})$. Then $\widehat{U} = \overline{9}(x_0 - \overline{3})^2$. Both the seed $\mathbf{s}$ given by the double root of $\widehat{U}$ and the seed $\mathbf{t}$ yield doubly infinite order-$5$ Somos sequences over $\mathbb{F}_{11}$ with period length $20$.

Once the analogues of Lemmas \ref{o6cm}--\ref{o6nz} are in place, the proofs of Theorems \ref{o4r}--\ref{o5t} can follow the same overall plan as before.

The material of Section \ref{int} admits lower-order analogues as well. We do not state them explicitly, as with orders $2$ and $3$ they are straightforward, and with orders $4$ and $5$ the derivations do not differ significantly from the ones in the setting of orders $6$ and $7$. For order $4$, we can apply Lemma \ref{df} with $k = 2$, $u_i = S_{i + 1}$, $e = (0, 2)$, $e' = (-2, 0)$, $e'' = (0, 2)$ as in the proof of Theorem \ref{o6i}. For order $5$, we can apply Lemma \ref{hdf} with $k = 4$, $u_i = S_{i + 2}$, $e = (0, 2)$, $e' = (-4, 0)$, $e'' = (0, 4)$ as in the proof of Theorem \ref{o7i}.

\section{Higher Orders and Further Work} \label{high}

The finite-rank properties of Somos sequences with $2 \le n \le 7$ are unlikely to generalise when $n \ge 8$. For example, direct computation over finite fields shows that the diamond rank of the order-$8$ unit Somos sequence and the half-diamond rank of the order-$9$ unit Somos sequence both exceed $2500$. (We use ``the diamond rank of $s$'' as shorthand for ``the diamond rank of $s \times s$'', and similarly with half-diamond ranks.)

So, are the order-$8$ unit and master Somos sequences of infinite diamond rank? What about the half-diamond ranks of the order-$9$ unit and master Somos sequences? Or, more generally, what about orders $n \ge 10$?

Still, there is one subclass of Somos sequences for which generalisations might hold after all. Let $\mathbf{n} = (n_1, n_2, n_3)$ with $n_1$, $n_2$, $n_3$ being positive integers such that $n = n_1 + n_2 + n_3$. We say that $s$ is a \emph{Gale-Robinson sequence} \cite{G} of order $n$ and type $\mathbf{n}$ when it satisfies \[s_is_{i + n} = a_1s_{i + n_1}s_{i + n_2 + n_3} + a_2s_{i + n_2}s_{i + n_3 + n_1} + a_3s_{i + n_3}s_{i + n_1 + n_2}\] for all $i$. So an order-$n$ Gale-Robinson sequence is a special case of an order-$n$ Somos sequence where almost all of the coefficients are set to zero.

Every Somos sequence of order $3 \le n \le 7$ is also a Gale-Robinson sequence, and much of the material in the first half of Section \ref{init} carries over to Gale-Robinson sequences. For example, every type $\mathbf{n}$ is associated with its own Gale-Robinson master sequence $S$.

The integrality properties discussed in Section \ref{int} are known \cite{G, FZ} to hold for all Gale-Robinson sequences. So it is natural to wonder if the finite-rank properties of low-order Somos sequences might generalise to arbitrary Gale-Robinson sequences as well.

Before we get to the experimental data, we must do some preliminary tidying up. Let $d = \gcd(n_1, n_2, n_3)$. If $d \ge 2$, we can split $s$ into $d$ decimations, by a factor of $d$, each one of which is a Gale-Robinson sequence itself, with the same coefficients but of type $\mathbf{n}/d$ instead. From now on, we will be interested most of all in the case when $d = 1$; i.e., when $n_1$, $n_2$, $n_3$ are relatively prime in aggregate. We call such types $\mathbf{n}$ \emph{primitive}.

The primitive-type condition of Conjectures \ref{gre} and \ref{gro} cannot be omitted. For example, with $\mathbf{n} = (2, 4, 6)$, $n = 12$, $\mathbf{a} = (1, 1, 1)$, $\mathbf{s} = (1, 1, \ldots, 1, 2)$, direct computation over finite fields shows that the diamond rank of $s$ exceeds $2500$; while $m = 5$ and $2^5 = 32$. Similarly, with $\mathbf{n} = (3, 6, 12)$, $n = 21$, $\mathbf{a} = (1, 1, 1)$, $\mathbf{s} = (1, 1, \ldots, 1, 2, 3)$, the half-diamond rank of $s$ exceeds $2500$ once again; while $m = 9$ and $2^9 = 512$.

Notice that for primitive types of order $n \ge 6$ we can assume without loss of generality that $n_1$, $n_2$, $n_3$ are pairwise distinct, in the sense that every Gale-Robinson sequence of a type which breaks this condition also belongs to another type which respects it. Explicitly, if $n_1 < n_2 = n_3$, set $\mathbf{n}^\star = (n_1, n_1 + n_2, n_2 - n_1)$; otherwise, if $n_1 > n_2 = n_3$, set $\mathbf{n}^\star = (2n_2, n_2, n_1 - n_2)$. Either way, if $s$ is a Gale-Robinson sequence of type $\mathbf{n}$ with coefficients $\mathbf{a} = (a_1, a_2, a_3)$, then it is also a Gale-Robinson sequence of type $\mathbf{n}^\star$ with coefficients $\mathbf{a}^\star = (a_1, a_2 + a_3, 0)$.

We call a type $\mathbf{n}$ \emph{proper} if it is primitive and $n_1$, $n_2$, $n_3$ are pairwise distinct. For Conjectures \ref{gre} and \ref{gro}, ``primitive'' and ``proper'' are interchangeable when $n \ge 6$. However, in other contexts it is sometimes more convenient to restrict consideration solely to the proper types.

We may approach Conjecture \ref{gre} experimentally as follows: Fix a type $\mathbf{n}$. Choose $\mathbf{a}$ and $\mathbf{s}$ uniformly at random out of, say, $[1; 100]^3$ and $[1; 100]^n$. Choose also a prime $p$ uniformly at random out of, say, all primes in the interval $[10^5; 10^6]$. Compute, next, many terms of the Gale-Robinson sequence $s$ over $\mathbb{F}_p$ determined by $\mathbf{a}$ and $\mathbf{s}$. Finally, take two large diamond sub-matrices in $s \times s$, of opposite parities, and compute their ranks over $\mathbb{F}_p$. The last step is nontrivial when these diamond sub-matrices are indeed very large; most of the experimental data reported in this section was obtained with the help of the FLINT software package \cite{FT}.

Suppose that, for many choices of $\mathbf{a}$, $\mathbf{s}$, and $p$, we consistently obtain one and the same rank $r$ which is much smaller than the sizes of the diamond sub-matrices being sampled. Then it would be reasonable to guess that $r$ ought to be the diamond rank of $S \times S$, for the corresponding Gale-Robinson master sequence $S$. Of course, Conjecture \ref{gro} can be approached analogously.

The author has run experiments of this kind for all proper types $\mathbf{n}$ with $8 \le n \le 25$, the experimental results being in full agreement with Conjectures \ref{gre} and \ref{gro}.

One subtlety is worth remarking upon. Define $m$ by $n = 2m + 2$ for even $n$ and $n = 2m + 3$ for odd $n$, as in Conjectures \ref{gre} and \ref{gro}. We refer to $2^m$ as the ``default'' rank of $\mathbf{n}$. We also refer to the value of $r$ indicated by the aforementioned series of experiments as the ``experimental'' rank of $\mathbf{n}$. For most proper types $\mathbf{n}$ with $8 \le n \le 25$, these two quantities coincide. However, in a few exceptional cases, the experimental rank is smaller.

Over the interval $8 \le n \le 25$, these exceptional types admit a simple description -- they are precisely the proper types $\mathbf{n}$ where two of $n_1$, $n_2$, $n_3$ share a greatest common divisor $g$ with $g \ge 5$. Furthermore, for all of these exceptional types, the ratio of the experimental rank to the default rank equals $\eta(g)$, where $\eta(2k + 1) = \eta(2k + 2) = (k + 1)/2^k$. This formula also agrees with the author's experimental results on some exceptional types $\mathbf{n}$ of higher orders $n \ge 26$. (Though it is difficult to predict how the formula might generalise to the proper types $\mathbf{n}$ for which two of the conditions $\gcd(n_1, n_2) \ge 5$, $\gcd(n_2, n_3) \ge 5$, $\gcd(n_3, n_1) \ge 5$ are satisfied simultaneously.)

So far, our discussion has been focused exclusively on the potential generalisations of Theorems \ref{o6r}~and~\ref{o7r}. However, some of the machinery we developed in order to prove these theorems might admit interesting generalisations as well.

We begin with the invariants. For any proper type $\mathbf{n}$, we can define the space $\mathcal{E}$ as in Section~\ref{init}; it is not too difficult to see that $\mathcal{E}$ will still depend only on the parity of $n$. Then, based on $\mathcal{E}$, we can define the subspace $\Upsilon_\boxtimes$ of $\Upsilon$ and the kernel $\Omega_\boxtimes$ of $\varphi$ over $\Upsilon_\boxtimes$ as in Section \ref{inv}. Notice that, in this setting, the distinction between primitive and proper types becomes meaningful.

\begin{conjecture} \label{cinv} For every proper type $\mathbf{n}$ of order $n$, it holds that $\dim \Omega_\boxtimes = \lfloor n/2 \rfloor$. \end{conjecture} 

So, for a Gale-Robinson sequence of a proper type $\mathbf{n}$, we would expect to find $m$ linearly independent nontrivial invariants. (Subtracting out the trivial invariant which corresponds to the element $\Pi$ of $\Omega_\boxtimes$.) We already know that Conjecture \ref{cinv} is true of all $3 \le n \le 7$. Direct computation confirms it also for both proper types of order $8$ as well as all three proper types of order $9$. The experimental data suggests the following supplements:

\begin{conjecture} \label{cinve} In the setting of Conjecture \ref{cinv}: (a) Every element $\Phi$ of $\Omega_\boxtimes$ exhibits the symmetry $\Phi(x_0, x_1, \ldots, x_{n - 1}) = \Phi(x_{n - 1}, x_{n - 2}, \ldots, x_0)$; and (b) There exists a basis of $\Omega_\boxtimes$ which consists entirely of positive-coefficient polynomials when viewed over $\mathbb{Z}[\boldsymbol{\alpha}, \mathbf{x}]$. \end{conjecture} 

It is well-known \cite{CS} that, in every Gale-Robinson master sequence, the numerators are positive-coefficient polynomials when viewed over $\mathbb{Z}[\boldsymbol{\alpha}, \mathbf{x}]$. This makes (b) somewhat more plausible. A different result of \cite{CS} establishes unit coefficients in certain three-dimensional arrays associated with the Gale-Robinson master sequences. It might be possible to strengthen Conjecture \ref{cinve} along similar lines, as outlined below.

For all proper types $\mathbf{n}$ with $3 \le n \le 9$, we can find a basis of $\Omega_\boxtimes$ where the coefficients of all basis polynomials are in the set $\{1, 2\}$ and the non-unit coefficients occur only at summands which are symmetric in the sense of (a). This suggests that it might be fruitful to express the elements of our basis in the form $\Phi = \Psi(x_0, x_1, \ldots, x_{n - 1}) + \Psi(x_{n - 1}, x_{n - 2}, \ldots, x_0)$. The natural strengthening of (b), in light of these observations, would be that we can always find a basis for $\Omega_\boxtimes$ where all of the $\Psi$'s are unit-coefficient polynomials when viewed over $\mathbb{Z}[\boldsymbol{\alpha}, \mathbf{x}]$.

We proceed now to relate the invariants to the finite-rank properties of Gale-Robinson sequences. In the context of a fixed proper type $\mathbf{n}$, we call two Gale-Robinson sequences \emph{twins} when they share the same coefficients and all invariants as in Conjecture \ref{cinv} agree over their coefficients and their seeds. Below, the term ``generically'' is used informally to mean ``generically for some reasonable notion of genericity''.

\begin{conjecture} \label{grte} Let $s$ and $t$ be two twinned Gale-Robinson sequences with nonzero terms, of a proper type and even order $n$ with $n = 2m + 2$. Then, generically, the matrix $s \times t$ is of diamond rank at most $2^m$. \end{conjecture}

\begin{conjecture} \label{grto} Let $s$ and $t$ be two twinned Gale-Robinson sequences with nonzero terms, of a proper type and odd order $n$ with $n = 2m + 3$. Then, generically, the matrix $s \times t$ is of half-diamond rank at most $2^m$. \end{conjecture}

Once again, the author has gathered experimental data over various finite fields in full agreement with Conjectures \ref{grte}~and~\ref{grto} for all proper types $\mathbf{n}$ of orders $8$ and $9$.

For each concrete proper type $\mathbf{n}$, we could in principle attempt to prove Conjectures \ref{gre}--\ref{grto} by carrying out computations similar to the ones in Sections \ref{o6} and \ref{o7}. However, already for orders $8$ and $9$, these computations become prohibitively difficult. Furthermore, such a strategy would at best allow us to handle only individual proper types $\mathbf{n}$ anyway; not all of them, or any infinite families of them. Clearly, deeper insights are necessary.

\section*{Acknowledgements}

The present paper was written in the course of the author's PhD studies under the supervision of Professor Imre Leader. The author is thankful to Prof.\ Leader for his unwavering support.

\end{document}